\documentclass[a4paper]{article}
\usepackage{latexsym, amsfonts, amsmath, amssymb}
\usepackage{tikz-cd}
\usepackage{microtype}
\usepackage[utf8]{inputenc}
\newcommand{\be}{\begin{equation}}
\newcommand{\ee}{\end{equation}}
\newcommand{\bea}{\begin{eqnarray}}
\newcommand{\eea}{\end{eqnarray}}
\newcommand{\bean}{\begin{eqnarray*}}
\newcommand{\eean}{\end{eqnarray*}}

\newcommand{\brray}{\begin{array}}
\newcommand{\erray}{\end{array}}
\newcommand{\ben}{\begin{equation}{nonumber}}
\newcommand{\een}{\end{equation}{nonumber}}

\newtheorem{dfn}{Definition}[section]
\newtheorem{thm}[dfn]{Theorem}
\newtheorem{lmma}[dfn]{Lemma}
\newtheorem{ppsn}[dfn]{Proposition}
\newtheorem{crlre}[dfn]{Corollary}
\newtheorem{xmpl}[dfn]{Example}
\newtheorem{rmrk}[dfn]{Remark}
\newcommand{\bdfn}{\begin{dfn}}
\newcommand{\bthm}{\begin{thm}}
\newcommand{\blmma}{\begin{lmma}}
\newcommand{\bppsn}{\begin{ppsn}}
\newcommand{\bcrlre}{\begin{crlre}}
\newcommand{\bxmpl}{\begin{xmpl}}
\newcommand{\brmrk}{\begin{rmrk}}
\newcommand{\edfn}{\end{dfn}}
\newcommand{\ethm}{\end{thm}}
\newcommand{\elmma}{\end{lmma}}
\newcommand{\eppsn}{\end{ppsn}}
\newcommand{\ecrlre}{\end{crlre}}
\newcommand{\exmpl}{\end{xmpl}}
\newcommand{\ermrk}{\end{rmrk}}
\newcommand{\IC}{\mathbb{C}}
\newcommand{\IR}{\mathbb{R}}

\newcommand{\oneformclassical}{\Omega^1 ( M )}

\newcommand{\twoform}{{\Omega}^2( \mathcal{A} )}
\newcommand{\tensora}{\otimes_{\mathcal{A}}}
\newcommand{\tensorsym}{\otimes^{{\rm sym}}_{\mathcal{A}}}
\newcommand{\tensorc}{\otimes_{\mathbb{C}}}
\newcommand{\A}{\mathcal{A}}
\newcommand{\C}{\mathcal{C}}

\newcommand{\E}{\mathcal{E}}
\newcommand{\F}{\mathcal{F}}
\newcommand{\G}{\mathcal{G}}
\newcommand{\Acenter}{\mathcal{Z}( \mathcal{A} )}
\newcommand{\Ecenter}{\mathcal{Z}( \mathcal{E} )}

\newcommand{\Psym}{P_{\rm sym}}
\newcommand{\Hom}{{\rm Hom}}
\newcommand{\id}{{\rm id}}

\newcommand{\clb}{{\cal B}}

\newcommand{\clh}{{\cal H}}

\def\a*{{\cal A}_{h,*}}
\def\B{{\cal B}(h)}
\def\B1{{\cal B}_1(h)}
\def\b{{\cal B}^{\rm s.a.}(h)}
\def\b1{{\cal B}^{\rm s.a.}_1(h)}

\def \qed {$\Box$}

\def\a*{{\cal A}_{h,*}}
\def\B{{\cal B}(h)}
\def\B1{{\cal B}_1(h)}
\def\b{{\cal B}^{\rm s.a.}(h)}
\def\b1{{\cal B}^{\rm s.a.}_1(h)}

\newcommand{\RNum}[1]{\uppercase\expandafter{\romannumeral #1\relax}}

\begin{document}
\begin{center}
{\Large{\bf A New Look at Levi-Civita connection in noncommutative geometry }}\\
\vspace{0.2in}
{\large {Jyotishman Bhowmick*, Debashish Goswami* and Soumalya Joardar**}}\\
*Indian Statistical Institute, Kolkata  \\
** IISER Kolkata, Mohanpur, Nadia, India \\
Emails: jyotishmanb$@$gmail.com, goswamid$@$isical.ac.in, soumalya.j$@$gmail.com \\
\end{center}
\begin{abstract}
 We prove the existence and uniqueness of Levi-Civita connections for strongly $\sigma$-compatible pseudo-Riemannian metrics on tame differential calculi. Such pseudo-Riemannian metrics properly contain the classes of bilinear metrics as well as their conformal deformations. This extends the previous results in \cite{article1} and \cite{article3}.  
\end{abstract}

\section{Introduction} 

We continue the study of Levi-Civita connection on noncommutative manifolds initiated in \cite{article1}, \cite{article3} and \cite{article4}. In these papers, the authors worked in the set up of tame differential calculus which allows one to define the notions of symmetry of a pseudo-Riemannian metric ( see Definiton \ref{metricdefn} ) and metric-compatibility of any connection ( see Definition \ref{19thfeb202021} ). Most importantly, we have the following theorem:
 \bthm ( \cite{article3}, \cite{article1} ) \label{9thaugust20}
 If $g$ is a pseudo-Riemannian {\bf bilinear} metric ( i.e., both left and right $\A$-linear ) on the bimodule of one-forms $ \Omega^1 ( \A ) $ of a tame differential calculus over an algebra $\A,$ then there exists a unique torsionless connection on $\Omega^1 ( \A ) $ which is also compatible with $g.$
\ethm

Thus, in \cite{article1}, \cite{article3}, \cite{article4},  our connections are defined on the level of one-forms $\Omega^1 ( \A ) $ of the differential calculus and a pseudo-Riemannian metric $g$ is a map from $ \Omega^1 ( \A ) \tensora \Omega^1 ( \A )  $ to $\A$ satisfying some properties. There is a parallel body of work where the connections are defined on suitable bimodules of derivations instead. This approach was pioneered by Rosenberg ( \cite{Rosenberg} ) followed by \cite{sheu}, \cite{pseudo} and more recently in \cite{tiger} and \cite{cylinder}. In \cite{article4}, it is proven that if $\Omega^1 ( \A ) $ denotes the bimodules of one-forms of a tame differential calculus, then there exists a canonical $\Acenter$-bimodule of derivations ( to be denoted by $\mathcal{X} ( \A ) $ ) contained in $\Omega^1 ( \A ) $ which plays the role of vector fields. In the presence of a {\bf bilinear} pseudo-Riemannnian metric on $\Omega^1 ( \A ), $ there is a one to one correspondence between connections on $ \Omega^1 ( \A ) $ and covariant derivatives on $ \mathcal{X} ( \A ). $ For more details, we refer to \cite{article4}.

On the other hand, there has been a lot of research activity around Levi-Civita connections on the level of forms. Bimodule connections were studied in great detail by Beggs, Majid and their collaborators for which re refer to \cite{beggsmajidbook}. In \cite{article3}, it has been shown that the Levi-Civita connection for a {\bf bilinear} pseudo-Riemannian metric on a tame differential calculus is actually a bimodule connection. Beggs and Majid  also studied ( \cite{majid_2} ) the issue of compatibility of the Levi-Civita connection in the presence of $\ast$-structures. For the question of existence of Levi-Civita connections on quantum groups and their homogeneous spaces, we refer to \cite{heckenberger_etal}, \cite{article6}, \cite{matassa},  \cite{landiqhom} and \cite{majid_1}. The existence and uniqueness of Levi-Civita connections for certain differential calculi over quasi-commutative algebras has been proven in \cite{aschieri} and \cite{weber}. Finally, for investigations in the case of finite spaces and groups, we refer to Chapter 8 of \cite{beggsmajidbook} and \cite{sitarz}.

The first new result of this article is the derivation of sufficient conditions for the Levi-Civita connections for {\bf bilinear} pseudo-Riemannian metrics on tame differential calculi are star-compatible and star-preserving in the sense of \cite{majid_2}. In particular, we prove that any bimodle connection on the space of one-forms on a tame differential calculus is automatically star-compatible. We also show that the Levi-Civita connections  for a natural pseudo-Riemannian bilinear metric on the noncommutative torus, quantum Heisenberg manifold and Cuntz algebras studied in \cite{article1} and \cite{soumalya} are star-preserving.

Our next goal is to generalize Theorem \ref{9thaugust20} for a larger class of pseudo-Riemannian metrics. In \cite{article3}, the proof of Theorem \ref{9thaugust20} was derived by imitating the classical proof of existence and uniqueness of Levi-Civita connections and yields a Koszul-type formula for this connection. We have been unable to generalize this proof to the case of pseudo-Riemannian metrics which are not necessarily bilinear. Nevertheless, the approach taken in \cite{article1} helps us to establish the main result of this article which is as follows:
\bthm ( Theorem \ref{existenceuniqueness} ) \label{9thaugust202}
 Suppose $ ( \Omega ( \A ), d ) $ is a tame differential calculus over an algebra $\A$ and $\Omega^1 ( \A ) $ be the bimodule of one-forms. If $g$ is a strongly $\sigma$-compatible pseudo-Riemannian metric on $ \Omega^1 ( \A ), $ then there exists a unique  connection on $ \Omega^1 ( \A ) $ which is torsionless and compatible with $g.$
 \ethm

We have imposed two different kinds of restriction in the above theorem. Let us discuss them one by one. Firstly, our proof works only in the set up of tame differential calculi ( see Definition \ref{tame} ). The isomorphism condition i. of Definition \ref{tame}  is restrictive. For example, it rules out  topologically nontrivial bimodules for algebras with trivial centers. However, we still have some interesting class of differential calculi which are tame. This is explained in Example \ref{exampleassumption}. 

The second restriction is on the variety of pseudo-Riemannian metrics which we allow. For the definition of strongly $\sigma$-compatible pseudo-Riemannian metrics, we refer to Definition \ref{stronglycompatible}. The class of such pseudo-Riemannian metrics properly include pseudo-Riemannian bilinear metrics as well as their conformal deformations ( see Proposition \ref{3rdjune202} ). While the bilinear case was resolved in \cite{article3} ( and \cite{article1} ), we have found a shorter proof for the case of conformal deformations of pseudo-Riemannnian bilinear metrics in \cite{article7}. However, let us point out that the technique of the proof of Theorem \ref{existenceuniqueness} inspired the proof of the main result in \cite{article6}. Moreover, since our main result gives a unified proof of both the bilinear case and the conformal deformation case, we believe that the proof of Theorem \ref{9thaugust202} merits an exposition.  
Following \cite{article3}, we have also defined the Ricci and scalar curvature of the Levi-Civita connection for any tame differential calculus in \cite{article7}. Thus, we have obtained the formulas of scalar curvature directly from the Levi-Civita connection as opposed to the approach of deriving the scalar curvature from the asymptotic expansion of the Laplacian in the seminal works  \cite{scalar_3}, \cite{Connes_moscovici}, \cite{khalkhali}  and references therein. 

	 We begin by recalling the definition and some properties of centered bimodules. In Section \ref{formulation}, we define the notion of tame differential calculus from \cite{article3}, \cite{article4} and then introduce pseudo-Riemannian metrics and metric compatibility of a connection. In Section \ref{starcomplc}, we recall the definitions of star-compatible and star-preserving connections from \cite{majid_2} and discuss the star-compatibility and star-preserving property of Levi-Civita connections studied in \cite{article1} and \cite{article3}.  Finally, we state and prove the main result giving the existence and uniqueness of the Levi-Civita connection for strongly $\sigma$-compatible pseudo-Riemannian metrics in Section \ref{exun}. 

We  fix some notations which we will follow. Throughout the article, $ \mathcal{A} $ will denote a complex algebra and $ \Acenter $ will denote its center. The tensor product over the complex numbers $ \IC $ is denoted by $ \tensorc $ while the notation $\tensora$ will denote the tensor product over the algebra $ \A. $    For a linear map $ T $ between suitable modules over $\A,$ $ {\rm Ran} ( T ) $ will denote the Range of $ T. $  

The class of bimodules over $\A$ with which we are interested will turn out to be centered bimodules ( see part i. of Lemma \ref{17thdec20192}  ). We recall the definition here.
\bdfn \label{centered}
We will say that a subset $ S  $ of a right $ \A $-module $ \E $ is right $ \A $-total in $ \E $ if   the right $\A$-linear span of $ S  $ equals $ \mathcal{E}. $ The center of  an $ \A $-bimodule $ \mathcal{E} $ is defined to be the set $ \mathcal{Z} ( \mathcal{E} ) = \{ e \in \mathcal{E}: e a = a e ~ \forall ~ a ~ \in \A  \}. $ It is easy to see that $ \mathcal{Z} ( \E ) $ is a $ \mathcal{Z} ( \A ) $-bimodule.  $ \mathcal{E} $ is called centered if $ \mathcal{Z} ( \mathcal{E} ) $ is right $ \A $-total in $ \E. $
\edfn
As an immediate corollary to the definition, we have the following simple lemma which we state without proof:
\blmma \label{centeredremark}
Suppose $ \E $ is a centered bimodule over $\A.$ Then the following statements hold:
  \begin{enumerate}
 \item[i.]   $ \Ecenter $ is also left $\A$-total in $\E.$

\item[ii.]  The set $ \{ \omega \tensora \eta: \omega, \eta \in \Ecenter \} $ is both left and right $\A$-total in $\E \tensora \E.$

\item[iii.]  If $X$ is an element of $\E \tensora \E,$ there exist $ v_i $ in $\E,$ $ w_i \in \Ecenter $ and $ a_i $ in $\A$ such that 
   $$ X = \sum_i v_i \tensora w_i a_i. $$
\end{enumerate} 	
\elmma

Let us demonstrate a simple example of a centered bimodule and compute its center.
\bxmpl \label{11thaugust20} 
 Suppose $( \E, d )$ is a tame differential calculus such that the bimodule $\E$ of one-forms is finitely generated and free as a right $\A$-module. Let $ \{ e_1, e_2, \cdots e_n \} $ be a basis of $\E$ such that $e_i \in \Ecenter $ for all $i.$ Then $\E$ is centered and 
$$ \Ecenter = \{  \sum_i e_i a_i: a_i \in \Acenter  \}. $$ 

Since $e_i$ belongs to $\Ecenter,$ it is clear that $\E$ is centered and in fact the set $ \{ \sum_i e_i a_i: a_i \in \Acenter \} \subseteq \Ecenter. $ Conversely, suppose that $ \sum_i e_i a_i $ is an element of $\Ecenter$ some elements $ a_1, a_2, \cdots a_n $ in $\A.$ Then for all $b$ in $\A,$ we get
$$ b ( \sum_i e_i a_i ) = ( \sum_i e_i a_i ) b $$
and so
$$ \sum_i e_i b a_i = \sum_i e_i a_i b. $$
Therefore, for all $i,$ the element $a_i$ commutes with all elements of $\A$ so that $a_i$ belongs to $\Acenter.$
\exmpl

For right $ \A  $ modules $ \mathcal{E} $ and $ \mathcal{F}, ~ {\rm Hom}_\A ( \mathcal{E}, \mathcal{F} ) $ will denote the set of all right $ \A $-linear maps from $ \E $ to $ \F. $ $ {\rm Hom}_\A ( \mathcal{E}, \mathcal{F} ) $ is an $\A$-bimodule in a natural way. The left $ \A $-module structure is given by left multiplication by elements of $ \A, $ i.e, for elements $ a $  in $ \A, ~ e  $ in $ \E $ and $ T  $ in $ {\rm Hom}_\A ( \mathcal{E}, \mathcal{F} ),$
 \begin{equation} \label{21stfeb20202} ( a. T ) ( e ) := a. T ( e ) \in \F. \end{equation} 
The right $ \A $ module structure on  $ {\rm Hom}_\A ( \mathcal{E}, \mathcal{F} ) $ is given by 
\begin{equation} \label{21stfeb20203} T. a ( e ) = T ( a e ). \end{equation}        
\bdfn \label{Estar}
If $ \E $ is an $\A$-bimodule, then $ \E^* $ will stand for the $\A$-bimodule $ {\rm Hom}_{\A} ( \mathcal{E}, \A ).$ 
 \edfn 
Finally, the following isomorphism will be used in the sequel multiple times and so we record it here.
\bdfn \label{xi}
Suppose $\E$ and $\F$ are finitely generated projective right $\A$-modules. The map $ \zeta_{\E, \F} $ will denote the canonical right $\A$-module isomorphism from $ \E \tensora \F^*  $ to $ \Hom_{\A} ( \F, \E ) $ which is defined by the following formula: 
$$ \zeta_{\E, \F} ( \sum_i e_i \tensora \phi_i ) ( f ) = \sum_i e_i \phi_i ( f ). $$
\edfn
The fact that $ \zeta_{\E, \F} $ is an isomorphism is well-known to the experts. We refer to Proposition 2.3 of \cite{article1} for a proof.

\section{Differential calculus and pseudo-Riemannian metrics} \label{formulation}
	
	In this section, we explain the set up of the article. The formulation of the Levi-Civita problem requires four ingredients: connections, metric, torsion of a connection and compatibility of a connection with a metric. We will introduce them one by one. Throughout this article, we will work with right connections on the bimodule of one-forms of a differential calculus ( see Definition \ref{connection} ) and the torsion is the same as the one studied by other authors. The novelties are in the definition of the pseudo-Riemannian metric and the metric-compatibility of a connection for which we have followed \cite{article1} and \cite{article3}.
	
	In order to make sense of the symmetry of the pseudo-Riemannian metric ( see \eqref{16thdec20191} below ), we need the notion of a quasi-tame differential calculus which is discussed in the first subsection. In the next subsection, we discuss tame differential calculi which we need for defining compatibility of a connection with a pseudo-Riemannian metric. It will turn out ( see Theorem \ref{torsionless} and Theorem \ref{existenceuniqueness} ) that a quasi-tame differential calculus guarantees the existence of a torsionless connection on the bimodule of one-forms while a tame differential calculus and the presence of a bilinear pseudo-Riemannian metric guarantees the existence of a unique Levi-Civita connection for a certain class of right $\A$-linear pseudo-Riemannian metric, which we call strongly $\sigma$-compatible ( see Definition \ref{stronglycompatible} ). The main result of \cite{article1} was the existence and uniqueness of a Levi-Civita connection for any {\bf bilinear} pseudo-Riemannian metric. The main result of this article ( Theorem \ref{existenceuniqueness} ) generalizes that result for strongly $\sigma$-compatible pseudo-Riemannian metrics. We end the section by recalling a necessary and sufficient condition for the existence and uniqueness of Levi-Civita connections. This is Theorem \ref{19thfeb2020}. Theorem \ref{existenceuniqueness} will be proved by verifying the hypothesis of Theorem \ref{19thfeb2020}.
	
	Let us recall the definition of a differential calculus.
	\bdfn \label{diffcal}
	Suppose $\A$ is an algebra over $\IC.$ A differential calculus on $\A$ is a triplet $ (  \Omega ( \A ), \wedge, d  ) $ such that the following conditions hold:
		 \begin{enumerate}
  	\item[i.] $ \Omega ( \A ) $ is an $\A$-bimodule,
			\item[ii.] $ \Omega ( \A ) = \oplus_{ i \geq 0 } \Omega^i ( \A ), $ where $ \Omega^0 ( \A ) = \A $ and $ \Omega^i ( \A ) $ are $\A$-bimodules. 
			\item[iii.] We have a bimodule map $ \wedge: \Omega ( \A ) \tensora \Omega ( \A ) \rightarrow \Omega ( \A ) $ such that 
			   $$ \wedge ( \Omega^i ( \A ) \tensora \Omega^j ( \A )  ) \subseteq \Omega^{i + j} ( \A ), $$
			\item[iv.] We have a map $ d: \Omega^i ( \A ) \rightarrow \Omega^{i + 1} ( \A ) $ such that 
				$$ d^2 = 0 ~{\rm and} ~ d ( \omega \wedge \eta ) = d \omega \wedge \eta + ( - 1 )^{{\rm deg} ( \omega ) } \omega \wedge d \eta, $$
			\item[v.] $ \Omega^i ( \A ) $ is spanned by elements of the form $ d a_1 \wedge \cdots \wedge d a_i a_{i + 1}. $
			\end{enumerate}
	  \edfn
	From now on, the bimodule of one-forms of a generic differential calculus will be denoted by $\E.$ We will always assume that $\E$ is finitely generated and projective as a right $\A$-module. For notational convenience, we will sometimes denote a differential calculus by a pair $ ( \E, d ) $ if $\E$ is the bimodule of one-forms of a differential calculus $ ( \Omega ( \A ), \wedge, d ). $

		Now we recall the definition of a right connection.
\bdfn ( \cite{frolich}, \cite{connes} ) \label{connection}
Suppose $ ( \Omega ( \A ), \wedge, d ) $ be a differential calculus on $\A.$  A (right) connection on $ \E:= \Omega^1 ( \A ) $
is a ${\mathbb C}$-linear map 
$$\nabla : \Omega^1 ( \A )   \rightarrow \Omega^1 ( \A )   \tensora \Omega^1 ( \A )  ~ {\rm satisfying} ~ {\rm the} ~ {\rm equation} ~  \nabla( e a) = \nabla( e )a + e \tensora da$$
for all $e$ in $\Omega^1 ( \A )   $ and $ a $ in $\A.$
 
The torsion of a connection $ \nabla $ is the right $\A$-linear map $ T_\nabla:= \wedge \nabla + d : \Omega^1 ( \A ) \rightarrow \Omega^2 ( \A ).$ The connection $\nabla$ is called torsionless if $ T_\nabla = 0. $
\edfn

\subsection{Quasi-tame differential calculi and pseudo-Riemannian metrics}

If $g$ is a pseudo-Riemannian metric on a smooth manifold $ M $ and $ \oneformclassical $ is the space of one-forms, then $g$ has to satisfy the equation:
\begin{equation} \label{16thdec20191} g ( e \tensora f ) = g ( f \tensora e ) \forall e, f \in \oneformclassical.  \end{equation} 
If $ {\rm  flip}: ~ \oneformclassical \otimes_{C^\infty ( M )} \oneformclassical \rightarrow \oneformclassical \otimes_{C^\infty ( M )} \oneformclassical $ denotes the map which sends $ e \otimes_{C^\infty ( M )} f $ to $ f \otimes_{C^\infty ( M )} e, $ then \eqref{16thdec20191} translates to the equation $ g \circ ~ {\rm flip} ~ = g. $ However, when we are dealing with a differential calculus over a noncommutative algebra $\A,$ then the map $ {\rm flip} $ is not well-defined. The framework of quasi-tame differential calculi gives rise to a bimodule map $\sigma$ which plays the role of the flip map in our set up. 
\bdfn \label{quasitame} ( Definition 2.11 of \cite{article1} )
Suppose $\E$ is the bimodule of one-forms of a differential calculus $ ( \Omega ( \A ), \wedge, d ).$ We say that the differential calculus is quasi-tame if the following conditions hold:
  \begin{enumerate}
  	\item[i.] The bimodule $\E$ is finitely generated and projective as a right $\A$ module.
  	\item[ii.] The following short exact sequence of right $\A$-modules splits:
		$$ 0 \rightarrow {\rm Ker} ( \wedge ) \rightarrow \E \tensora \E \rightarrow \Omega^2 ( \A ) \rightarrow 0. $$
		Thus, in particular, there exists a right $\A$-module $\F$ isomorphic to $\twoform$ such that:
  \begin{equation} \label{splitting25thmay2018}
  \E \tensora \E = {\rm Ker}(\wedge) \oplus \F 
  \end{equation}
    \item[iii.] Let us denote the idempotent in $ {\rm Hom}_{\A}(\E \tensora \E, \E \tensora \E)$ with range  ${\rm Ker}(\wedge)$ and  kernel $\F$ by the symbol $ \Psym. $ We assume that $ \Psym $ is an $\A$-bimodule map.
\end{enumerate}
 \edfn
Now we are ready to introduce the map $ \sigma $ which will play the role of the flip map in this article.  
\bdfn \label{defnpsymandsigma}
 If $ ( \E, d ) $ is a quasi-tame differential calculus, $ \sigma $ will denote the map $ 2  \Psym - 1. $
\edfn
Using the map $\sigma,$ we can now define pseudo-Riemannian metrics on a quasi-tame differential calculus. We will need the notation $ \E^*:= \Hom_\A ( \E, \A ) $ introduced in Definition \ref{Estar}.
\bdfn \label{metricdefn} ( \cite{article1}, \cite{article3} )
Suppose $\E$ is the bimodule of one-forms of a quasi-tame differential calculus $ ( \Omega ( \A ), \wedge, d ).$ A pseudo-Riemannian metric $ g $ on $ \E $ is
 an element of $ {\rm Hom}_{\A} ( \E \tensora \E, \A ) $ such that
 \begin{enumerate}
 \item[i.] $g$ is symmetric, i.e. $ g \circ \sigma = g, $
  \item[ii.] $g$  is non-degenerate, i.e,  the right $ \A$-linear map $ V_g: \E \rightarrow {\E}^* $ defined by $ V_g ( \omega ) ( \eta ) = g ( \omega \otimes_{\A} \eta ) $ is
 an isomorphism of right $ \A$-modules.
\end{enumerate}
We will say that a pseudo-Riemannian metric $g$ is a pseudo-Riemannian {\bf bilinear} metric if in addition, $ g $ is also left $\A$-linear.
\edfn
 Here, the right $ \A $-module structure on $ \E^* = \Hom_\A ( \E, \A ) $ is as in \eqref{21stfeb20203}, i.e, if $ \phi $ belongs to  $ \E^*, $ then for all $ a $ in  $ \A $ and $ e  $ in $\E,$ $ ( \phi. a ) ( e ) = \phi ( a. e ).  $

Since the map $V_g$ is one one, we have the following observation:
\brmrk \label{15thdec20191}
The condition ii. of Definition \ref{metricdefn} implies that if $g$ is a pseudo-Riemannian metric on a quasi-tame differential calculus and $e$ is an element of $\E$ such that $ g ( e \tensora f ) = 0 $ for all $f$ in $\E,$ then $e = 0.$ 
\ermrk
Since we are concerned about existence of Levi-Civita connection, it is important to know whether there exists at least one torsionless connection on $\E.$ The following theorem proved in \cite{article3} ( also see \cite{article1} ) answers the question in the affirmative. 
\bthm  \label{torsionless} ( Theorem 3.3 of \cite{article3} )
Suppose $ ( \Omega ( \A ), \wedge, d ) $ is a quasi-tame differential calculus. Then the bimodule of one-forms $ \E = \Omega^1 ( \A ) $ admits a torsionless connection.
\ethm

\subsection{Tame differential calculi and Levi-Civita connections} \label{metriccompatibility}

In this subsection, we discuss the notion of tame differential calculus and its interactions with pseudo-Riemannian metrics and connections.
Before defining a tame differential calculus, let us recall that if $ ( \E, d ) $ is a differential calculus such that the condition ii. of Definition \ref{quasitame} holds, then $ \E \tensora \E = {\rm Ker}(\wedge) \oplus \F $ for some right $\A$-module $\F$ and  we have an idempotent $\Psym$ in $ \Hom_\A ( \E \tensora \E, \E \tensora \E ) $ with range $ {\rm Ker} ( \wedge ) $ and kernel equal to $\F.$ The map $\sigma = 2 \Psym - 1 $ defined in Definition \ref{defnpsymandsigma} plays the role of the flip map.  
\bdfn \label{tame}
A differential calculus $ ( \E, d ) $ is called tame if the conditions {\rm i.} and {\rm ii.} of Definition \ref{quasitame} hold and moreover:
	\begin{enumerate}
		\item[i.] The map $u^\E:\Ecenter \otimes_{\Acenter} \A \rightarrow \E$ defined by
		$$u^\E(\sum_i e_i \otimes_{\Acenter} a_i)=\sum_i e_i a_i$$
		is an isomorphism of vector spaces,
		\item[ii.] $\sigma $ satisfies the following equation for all $\omega, \eta \in \Ecenter:$
		 \begin{equation} \label{17thdec20191} \sigma ( \omega \tensora \eta ) = \eta \tensora \omega. \end{equation}
	\end{enumerate}
\edfn
Let us note that a tame differential calculus is a quasi-tame differential calculus.
\brmrk \label{17thdec2019remark}
In Lemma 4.4 of \cite{article1}, it has been proved that the map $ \Psym $ is $\A$-bilinear and so a tame differential calculus is  quasi-tame. The proof of Lemma 4.5 of \cite{article1} is of course written in the set up of spectral triples but the proof goes verbatim for a differential calculus.
\ermrk

Before discussing the consequences for a differential calculus to be tame, we list a class of examples of differential calculi which are tame.

	\bxmpl \label{exampleassumption}
	
	\begin{enumerate}
	
	\item Let $ M $ be a Riemannian manifold. Consider the usual differential calculus $ ( \Omega^{\cdot} ( M ), \wedge, d   ) $ where $ \Omega^{\cdot} ( M ) = \oplus_i \Omega^i ( M ), $ $ \wedge $ is the usual wedge map and $ d $ the de-Rham differential. Let $ \A = C^\infty ( M ). $ Then $ \Omega^{\cdot} ( M ) $ is an $\A$-bimodule and  $ \E = \oneformclassical. $ Therefore, $ \mathcal{Z} ( \E ) = \oneformclassical, ~ \mathcal{Z} ( \A ) = \A = C^\infty(M). $ Here, $ \sigma $ is the usual flip map and all the conditions of Definition \ref{tame} are satisfied.

	\item Consider the differential calculi for the fuzzy $3$-sphere as in \cite{frolich},  quantum Heisenberg manifold ( \cite{rieffel_heisenberg}, \cite{chak_sinha} )  and the Rieffel deformations ( \cite{rieffel} ) of a classical compact Riemannian manifold by an isometric and free toral action as discussed in \cite{Connes-dubois} and \cite{connes_landi}. These differential calculi are also tame. We refer to Theorem 5.4,  Theorem 6.6 and Theorem 7.1 of \cite{article1} for the proofs.  Indeed, the proof of these results contain the fact that the bimodule of one forms for each of these examples satisfy the hypothesis of Theorem 4.1 of \cite{article1}. However, the hypotheses of Theorem 4.1 of \cite{article1} are exactly the conditions which make a differential calculus tame.
		
\item The example of the fuzzy $2$-sphere considered in \cite{article3} is another example of a tame differential calculus. We refer to Theorem 8.5 of \cite{article3} for the proof.

\item  Theorem 3.4 of \cite{soumalya} proves that a differential calculus constructed on the Cuntz algebra ( from a natural $C^*$-dynamical system )  is also tame.
\end{enumerate}	
			\exmpl
			
Let us now state some properties of a tame differential calculus. Lemma \ref{17thdec20192} derives some useful formulas of the maps $ \Psym $ and $ \sigma $ which will be used repeatedly in the sequel. In particular, \eqref{10thjuly20182} shows that $ \sigma $ is an analogue of the flip map. Lemma \ref{lemma0} states some properties of a pseudo-Riemannian metric on a tame differential calculus. 			

The following lemma needs the notion of centered bimodule which was defined in Definition \ref{centered}.
\blmma \label{17thdec20192}
Suppose $( \E, d )$ is a tame differential calculus. Then  the following statements hold:
\begin{enumerate}
		\item[i.] The bimodule $\E$ is centered.
		\item[ii.] For all $\omega \in \Ecenter$ and $e \in \E,$ we have
		    \begin{equation} \label{10thjuly20182} \sigma ( \omega \tensora e ) = e \tensora \omega, ~ \sigma ( e \tensora \omega ) = \omega \tensora e. \end{equation}
		\item[iii.] If $\Psym$ denotes the map defined in Definition \ref{quasitame}, then for all $\omega \in \Ecenter$ and $e \in \E,$
		   \begin{equation} \label{21stmarch20} \Psym ( e \tensora \omega ) = \Psym ( \omega \tensora e ) = \frac{1}{2} ( \omega \tensora e + e \tensora \omega ).\end{equation}
			In fact, the decomposition $ \E \tensora \E = {\rm Ker} ( \wedge ) \oplus \mathcal{F} $ on simple tensors is explicitly given by
			\begin{equation} \label{21stmarch202} \omega \tensora \eta a = \frac{1}{2} ( \omega \tensora \eta a + \eta \tensora \omega a  ) + \frac{1}{2} ( \omega \tensora \eta a - \eta \tensora \omega a  ) \end{equation}
			for all $\omega, \eta$ in $\Ecenter$ and for all $a$ in $\A.$
			\item[iv.] The map $\sigma$ is both left and right $\A$-linear.
\end{enumerate}		
\elmma
{\bf Proof:} Since the map $ u^{\E}: \Ecenter \otimes_{\Acenter} \A \rightarrow \E $ is an isomorphism of vector spaces, Proposition 2.4 of \cite{article1} implies that $\E$ is centered. This proves the first assertion. Now we prove the second assertion. Since $\E$ is centered by i., there exist elements $\eta_i \in \Ecenter $ and $ a_i $ in $\A$ such that $ e = \sum_i \eta_i a_i. $ Hence, by the right $\A$-linearity of $\sigma,$ we obtain
\begin{eqnarray*}
\sigma ( \omega \tensora e ) &=& \sum_i \sigma ( \omega \tensora \eta_i ) a_i = \sum_i \eta_i \tensora \omega a_i ~  {\rm (} ~ {\rm by} ~ \eqref{17thdec20191} ~ {\rm )}\\
                             &=& \sum_i \eta_i \tensora a_i \omega = \sum_i \eta_i a_i \tensora \omega = e \tensora \omega   
\end{eqnarray*}
since $ \omega a_i = a_i \omega $ as $\omega$ belongs to $\Ecenter.$ The second equation of \eqref{10thjuly20182} can be proved similarly.

Now we prove the third assertion. The equation \eqref{21stmarch20} follows from \eqref{10thjuly20182} since $ \Psym = \frac{1 + \sigma}{2}.$ On the other hand, the equation \eqref{21stmarch202} is already proved in part i. of Proposition 4.3 of \cite{article1}.

Finally, for iv., a combination of Remark \ref{17thdec2019remark} and the equality $ \sigma = 2 \Psym - 1 $ proves that $\sigma$ is left $\A$-linear. $\sigma$ is right $\A$-linear by the definition of $\Psym$ in Definition \ref{quasitame}.
\qed

As a simple consequence of the fact that $\E$ is centered, we have the following fact which was proved in Section 4 of \cite{article3} ( see equation (2) ):
\begin{equation} \label{18thdec20191} a. e = e. a ~ {\rm if} ~ e \in \E ~ {\rm and} ~ a \in \Acenter.  \end{equation}
This says that $\E$ is a central bimodule in the sense of \cite{dubois} and \cite{dubois2}.

Now we come to pseudo-Riemannian metrics on a tame differential calculus. Firstly, by Remark \ref{17thdec2019remark}, a tame differential calculus is quasi-tame and so the notion of a pseudo-Riemannian metric makes sense on a tame differential calculus. 

\blmma \label{lemma0}
 Suppose $ ( \E, d ) $ is a tame differential calculus and $g$ is a pseudo-Riemannian metric on $\E.$ Then the following are true:
\begin{itemize}
\item[i.] If  either $ e $  or $ f $ belongs to $ \Ecenter,$ we have
\begin{equation} \label{gsigmaisg} g (  e \tensora f ) = g ( f \tensora e ). \end{equation}
\item[ii.] If $ g_0 $ is a  pseudo-Riemannian  bilinear metric, then $ g_0 ( \omega \tensora \eta ) \in \mathcal{Z} ( \A ) $ if $ \omega, \eta  $ belong to $ \mathcal{Z} ( \E ).$
\item[iii.] If $e$ is an element of $\E$ such that $ g ( e \tensora \omega ) = 0 $ for all $ \omega  $ in $\Ecenter,$ then $e = 0.$ The same conclusion holds if $ g ( \omega \tensora e ) = 0 $ for all $ \omega $ in $\Ecenter.$
\item[iv.] $V_g$ is   left $\mathcal{Z} ( \A )$-linear.
\end{itemize}
\elmma
{\bf Proof:} The first assertion follows as a trivial consequence of \eqref{10thjuly20182} and the relation $ g \sigma = g. $ 
The second assertion was already proved in Lemma 3.6 of \cite{article1}. 
	
	Now we prove iii. Since $\E$ is centered ( Lemma \ref{17thdec20192} ), $\Ecenter$ is right $\A$-total in $\E.$ Therefore, an element $f$ in $\E$ can be written as a finite linear combination $ f = \sum_i \omega_i a_i $ for elements $ \omega_i $ in $\Ecenter$ and $ a_i  $ in $\A.$ Hence, if $ g ( e \tensora \omega ) = 0 $ for all $ \omega $ in $\Ecenter,$ then
	$$ g ( e \tensora f ) = \sum_i g ( e \tensora \omega_i ) a_i = 0 $$
	where we have used the right $\A$-linearity of $g.$ Thus, $ e = 0 $ by Remark \ref{15thdec20191}.
	
	Now suppose that  $ g ( \omega \tensora e ) = 0 $ for all $\omega$ in $\Ecenter.$ But \eqref{gsigmaisg} implies that $ g ( e \tensora \omega ) = 0 $ and we are back to the previous case.

 Finally we prove iv. We will use the left $\A$-module structure on $\E^* = \Hom_\A ( \E, \A ) $ as in \eqref{21stfeb20202}, i.e, for $ \phi $ in $\E^*,$ $a$ in $\A$ and $e$ in $\E,$ $ ( a. \phi ) ( e ) = a. \phi ( e ). $

  Let  $a^\prime \in \mathcal{Z} ( \A ),$ $e, f \in \E.$ We write $ f = \sum_i \omega_i a_i $ as in the proof of part iii. and use \eqref{18thdec20191} repeatedly  to compute:
 \begin{eqnarray*}
  V_g(a^\prime e)( f )&=& g(a^\prime e \tensora f )=g( e a^\prime \tensora \sum_i ( \omega_i a_i ) )= \sum_i g( e \tensora a^\prime \omega_i) a_i\\
                      &=& \sum_i g( e \tensora \omega_i a^\prime) a_i = \sum_i g( e \tensora \omega_i ) a_i a^\prime =  g( e \tensora f) a^\prime\\
											 &=& a^\prime g ( e \tensora f ) = a^\prime V_g( e )( f )
 \end{eqnarray*}
 where we have used the right $\A$-linearity of $g$ and  the fact that $a^\prime \in  {\mathcal Z}(\A)$. This proves iv. and finishes the proof of the lemma. \qed

The following example of a pseudo-Riemannian metric ( which is not necessarily bilinear )  will be used in Proposition \ref{3rdjune202}.

\bppsn \label{11thaugust202}
Suppose $( \E, d )$ is a tame differential calculus such that the bimodule $\E$ of one-forms is finitely generated and free as a right $\A$-module. Let $ \{ e_1, e_2, \cdots e_n \} $ be a basis of $\E$ such that $e_i \in \Ecenter $ for all $i.$ If $k$ is an invertible element of $\A,$ then the map 
$$ g:  \E \tensora \E \rightarrow \A, g ( ( \sum_i e_i a_i  ) \tensora ( \sum_j e_j b_j  ) ) = \sum_i k^i a_i b_i  $$
is a pseudo-Riemannian metric on $\E.$
\eppsn
{\bf Proof:} Let us note that  $ g ( e_i \tensora e_j ) = k^i \delta_{ij.}$ It is clear that $g$ is right $\A$-linear. 
In order to prove that $g$ is a pseudo-Riemannian metric, we need to check that $ g \circ \sigma = g $ and $ V_{g}: \E \rightarrow \E^* $ is a right $\A$-linear isomorphism ( see Definition \ref{metricdefn} ).
		We compute
		\begin{eqnarray*}
		g \circ \sigma (  ( \sum_i e_i a_i   ) \tensora ( \sum_j e_j b_j ) ) &=& \sum_{i,j} g \circ \sigma ( e_i \tensora e_j  ) a_i b_j = \sum_{i,j} g ( e_j \tensora e_i ) a_i b_j ~ {\rm (} ~ {\rm by} ~ \eqref{17thdec20191} ~ {\rm )}\\
	&=& \sum_{i,j} k^i \delta_{ij} a_i b_j = \sum_i k^i a_i b_i \\
	&=& g (  ( \sum_i e_i a_i   ) \tensora ( \sum_j e_j b_j ) ).
	\end{eqnarray*}
	Next, in order to prove that the map $ V_{g}: \E \rightarrow \E^* $ is an isomorphism, we note that $\E^*$ is a free right $\A$-module of rank $n$ generated by $ \{ \phi_1, \cdots \phi_n \} $ where $ \phi_j $ belonging to $\E^*$ is defined by
		$$ \phi_j ( \sum_i e_i a_i ) = a_j. $$
		We note that for all $a, b_i$ in $\A,$
		$$ ( \phi_j. a ) ( \sum_i e_i b_i ) = \sum_i \phi_j ( a e_i b_i ) = \sum_i \phi_j ( e_i a b_i ) = a b_j = ( a \phi_j ) ( \sum_i e_i b_i )  $$
		proving that $ \phi_j $ belongs to $\mathcal{Z} ( \E^* ). $
		
	Let us recall ( Definition \ref{metricdefn} ) that $V_{g}: \E \rightarrow \E^* $ is defined by $ V_{g} ( e ) ( f ) = g ( e \tensora f ) $ for all $e, f $ in $\E.$ Since $ g ( e_i \tensora e_j ) = k^i \delta_{ij}, $ it immediately follows that 
	\begin{equation} \label{11thaugust203} V_{g} ( e_i ) = k^i \phi_i = \phi_i k^i \end{equation}
	as $\phi_i$ belongs to $ \mathcal{Z} ( \E^* ). $ As $ V_{g} $ is right $\A$-linear and $k$ is invertible, it follows from $ V_{g} $ is onto. 
	
	For proving that $ V_{g} $ is one-one, let us assume that there exist elements $ a_1, \cdots  a_n $ in $\A$ such that $ V_{g} ( \sum_i  e_i a_i ) = 0. $ Therefore, for all $j,$
	$$ V_g ( \sum_i e_i a_i ) ( e_j ) = 0. $$
	But by \eqref{11thaugust203} and the fact that $e_i \in \Ecenter, $ we obtain
	$$ V_g ( \sum_i e_i a_i ) ( e_j ) = \sum_i \phi_i k^i a_i ( e_j  ) = \sum_i \phi_i ( e_j ) k^i a_i = k^j a_j. $$
	Since $k$ is invertible, this proves that $a_j = 0$ for all $j$ and so $V_g$ is one-one.
\qed

\brmrk \label{18thjan21}
Suppose $(\E, d )$ be a differential calculus as in the statement of Proposition \ref{11thaugust202}. If we put $k = 1$ in that proposition, then the pseudo-Riemannian metric $g$ is actually $\A$-bilinear.
\ermrk

 Let $ ( \E, d ) $ be a tame differential calculus, $g$ a pseudo-Riemannian metric on $\E$ and $\nabla$ a connection on $\E.$ We can now define the notion of compatibility of $\nabla$ with $g.$  As observed in Remark 4.15 of \cite{article1}, the proof of the following proposition works for any pseudo-Riemannian metric ( which is {\bf not necessarily bilinear} ).

\bppsn ( Subsection 4.1, \cite{article1} ) \label{30thmarch203}
For a pseudo-Riemannian metric on $\E,$ let us  define $\Pi_g^0(\nabla):\Ecenter\tensorc\Ecenter \rightarrow \E$ as the map given by
$$\Pi_g^0(\nabla)(\omega \tensorc \eta)=(g\tensora {\rm id})\sigma_{23}(\nabla(\omega)\tensora \eta + \nabla(\eta) \tensora \omega).$$
Then $ \Pi_g^0 $ extends to a well defined map from $ \E \tensora \E $ to $ \E $ to be denoted by $ \Pi_g ( \nabla ).$ 
\eppsn
In fact, if $\omega, \eta $ belong to $\Ecenter$ and $ a $ belongs to $\A,$ then 
 \begin{equation} \label{20thfeb2020} \Pi_g ( \nabla ) ( \omega \tensora \eta a ) = \Pi^0_g ( \nabla ) ( \omega \otimes_{\Acenter} \eta ) a + g ( \omega \tensora \eta ) da. \end{equation}
The equation \eqref{20thfeb2020} defines the map $ \Pi_g ( \nabla ) $ on the whole of $ \E \tensora \E.$ Indeed, by Lemma \ref{17thdec20192}, $\E$ is centered and so by Lemma \ref{centeredremark}, any element of $\E \tensora \E $ is a finite sum of elements of the form $ \omega \tensora \eta a, $ where $\omega, \eta \in \Ecenter$ and $a$ belongs to $\A.$  

\bdfn \label{19thfeb202021}
Suppose $ ( \E, d ) $ is a tame differential calculus and $ g $ is a pseudo-Riemannian metric on $\E.$

 A connection $\nabla$ on $\E$ is said to be compatible with $g$ if 
$$\Pi_g(\nabla) ( e \tensora f ) = d ( g ( e \tensora f ) )~ {\rm for} ~ {\rm all} ~ e,~ f ~ {\rm in} ~ \E.$$
A connection $\nabla$ on $\E$ which is torsionless and compatible with $g$ is called a Levi-Civita connection for the triplet $(\E, d, g ).$
\edfn

We end the section by stating two results from \cite{article1} and \cite{article3}. The first is a necessary and sufficient condition for the existence and uniqueness of Levi-Civita connection on a tame differential calculus and will be used to prove our main Theorem \ref{existenceuniqueness}. We will need a couple of definitions. 

\bdfn \label{19thfeb20202}
For a tame differential calculus $ ( \E, d ), $ the symbol $ \E \tensorsym \E$ will denote   ${\rm Ker} ( \wedge ) = {\rm Ran} ( \Psym ). $ If $g$ is a pseudo-Riemannian metric, the element $dg$ will denote the map
$$ dg: \E \tensora \E \rightarrow \E, dg ( e \tensora f ) = d ( g ( e \tensora f )  ).$$
\edfn
Let us note the following fact.
\brmrk \label{21stfeb20204}
From Definition \ref{diffcal}, the map $ \wedge: \E \tensora \E \rightarrow \twoform $ is an $\A$-bimodule map and so $ \E \tensorsym \E =  {\rm Ker} ( \wedge ) $ is an $\A$-bimodule.
\ermrk
We have the following theorem:
\bthm ( Theorem 4.14, \cite{article1} ) \label{19thfeb2020}
 Let $ ( \E, d ) $ be a tame differential calculus and $g$ a pseudo-Riemannian metric on $\E.$ Let $ \E \tensorsym \E $ be the $\A$-bimodule of Definition \ref{19thfeb20202}.  
$$ {\rm We} ~ {\rm define} ~ {\rm a} ~ {\rm map} ~ \Phi_g : {\rm Hom}_{\A}(\E, \E \tensorsym \E) \rightarrow {\rm Hom}_{\A}(\E \tensorsym \E, \E) ~ {\rm by} ~ {\rm the} ~ {\rm formula:} $$
$$ \Phi_g(L) ( X ) =(g \tensora {\rm id}) \sigma_{23} (L \tensora {\rm id})(1+\sigma) ( X ) $$
for all $X$ in $\E \tensorsym \E.$

Then $ \Phi_g $ is right $\A$-linear. Moreover, if $ \Phi_g:{\rm Hom}_{\A}(\E, \E \tensorsym \E) \rightarrow {\rm Hom}_{\A}(\E \tensorsym \E, \E) $ is an isomorphism of right $\A$-modules, then there exists a unique connection $ \nabla $ on $\E$ which is torsion-less and compatible with $g.$
Moreover, if $\nabla_0$ is a fixed torsionless connection on $\E,$ then  $ \nabla $ is given by the following equation: 
				\be \label{26thnov2} \nabla = \nabla_0 + \Phi^{- 1}_g ( dg - \Pi_g ( \nabla_0 ) ).  \ee
	Here, $dg: \E \tensora \E \rightarrow \E $ is the map defined in Definition \ref{19thfeb20202}.			
\ethm

The proof of this theorem works for any pseudo-Riemannian metric. The formula \eqref{26thnov2} follows from the proof of Theorem 4.14 of \cite{article1}. We only need to remark that the proof of Theorem 4.13 of \cite{article1} uses the existence of a torsion-less connection on $\E.$ In our case, this condition is satisfied by virtue of Theorem \ref{torsionless} since our differential calculus is in particular quasi-tame ( Remark \ref{17thdec2019remark} ). 

The main result of \cite{article3} and \cite{article1} is the following:

\bthm  ( Theorem 6.1 of \cite{article3}, Theorem 4.1 of \cite{article1} ) \label{21stfeb2020}
Let $ ( \E, d ) $ be a tame differential calculus and $g_0$ be a pseudo-Riemannian {\bf bilinear} metric on $\E.$ Then there exists a unique Levi-Civita connection for the triplet $(\E, d, g_0 ).$
\ethm
\brmrk \label{22ndfeb2020}
In \cite{article1}, Theorem \ref{21stfeb2020} was proved by verifying that the map $ \Phi_{g_0}: {\rm Hom}_{\A}(\E, \E \tensorsym \E) \rightarrow {\rm Hom}_{\A}(\E \tensorsym \E, \E) $ is an isomorphism of right $\A$-modules.
\ermrk
 In \cite{article3}, a completely different proof was given. Indeed, the uniqueness of such a connection followed by deriving a Koszul type formula of a torsionless and $g_0$ compatible connection. The existence followed by proving that the above mentioned Koszul-formula indeed defines a torsionless and $g_0$-compatible connection on $\E.$

For proving the existence and uniqueness of Levi-Civita connection ( on a tame differential calculus ) for  strongly $\sigma$-compatible pseudo-Riemannian metrics ( see Definition \ref{stronglycompatible} )  we will use Theorem \ref{19thfeb2020}. We have been unable to generalize the proof of \cite{article3} for such metrics. 

\section{Star-compatibility of Levi-Civita connections} \label{starcomplc}

In this section, we discuss the issue of star-compatibility of the Levi-Civita connection on a tame differential calculus. We will need some terminologies from Section 3 of \cite{majid_2}. However, since we are working with right connections ( as opposed to left connections, as in \cite{majid_2} ), our formulas will be slightly different. This is the only section in this article where we use star-structures.

 We will begin with the definition of bimodule connections.
\bdfn
Suppose $ ( \E, d ) $ is a differential calculus and $\F$ is an $\A$-bimodule.  
\begin{enumerate} 
\item Suppose $ \nabla $ is a right connection on $\F$ and $ \sigma^\prime : \E \tensora \F \rightarrow \F \tensora \E $ is an invertible $\A$-bilinear map.  Then the pair $ ( \nabla, \sigma^\prime )  $ is called a right bimodule  connection if
  $$ \nabla ( a f ) = a \nabla ( f ) + \sigma ( da \tensora f ) ~ {\rm for} ~ {\rm all} ~ f \in \F ~ {\rm and} ~ a \in \A. $$
\item Suppose $\nabla: \F \rightarrow \E \tensora \F $ is a  left connection on $ \F, $  i.e,  
$ \nabla ( a f ) = a \nabla ( f ) + da \tensora f $   for  all $ f \in \F $  and $ a \in \A. $ Moreover, let $ \sigma^\prime : \F \tensora \E \rightarrow \E \tensora \F $ be an invertible $\A$-bilinear map. 	Then the pair $ ( \nabla, \sigma^\prime ) $ is called a left bimodule connection if
$$ \nabla ( f a ) = \nabla ( f ) a + \sigma^\prime ( f \tensora da ). $$ 	
\end{enumerate} 
\edfn
The statements in the following proposition will be useful for us:
\bppsn ( Lemma 3.70 and Theorem 3.78, \cite{beggsmajidbook} ) \label{27thsep20}
Suppose $\F$ and $\G$ are $\A$-bimodules and $(\E, d )$ is a differential calculus on $\A.$
\begin{enumerate}
\item If $ ( \nabla, \sigma^\prime ) $ is a left bimodule connection on $\F,$ then $ ( ( \sigma^\prime )^{-1} \nabla, ( \sigma^\prime )^{-1} ) $ is a right bimodule connection.

\item If $ ( \nabla_1, \sigma_1 ) $ and $ ( \nabla_2, \sigma_2 ) $ are right bimodule connections on $\F$ and $\G$ respectively and
  $$ \nabla: ( \F \tensora \G ) \rightarrow ( \F \tensora \G ) \tensora \E, ~ \nabla ( f \tensora g ) = ( \sigma_2 )_{23} ( \nabla_1 ( f ) \tensora g ) + f \tensora \nabla_2 ( g ), $$
	then $ ( \nabla, ( \id \tensora \sigma_2 ) ( \sigma_1 \tensora \id ) ) $ is a right bimodule connection.	
\end{enumerate}
\eppsn

Now we explain the formulation of star-compatibility of a bimodule connection as in \cite{majid_2}. Throughout this section, if $\A$ is a $\ast$-algebra and $\F$ is an $\A$-bimodule, then $\overline{\F}$ will denote the $\A$-bimodule whose underlying set is equal to $\F$ but equipped with a different $\A$-bimodule structure:
$$ a \overline{f} = \overline{f a^*}, ~ \overline{f} a = \overline{a^* f} ~ {\rm for} ~ {\rm all} ~ a \in \A, ~ f \in \F. $$
We will let $ ( {}_{\A} {\bf \mathcal{M}} {}_{\A}, \otimes ) $ denote the monoidal category of $\A$-bimodules, i.e, objects of $ {}_{\A} {\bf \mathcal{M}} {}_{\A} $ are $\A$-bimodules while morphisms are $\A$-bimodule maps. Let us make a list of  morphisms in the category $ {}_{\A} {\bf \mathcal{M}} {}_{\A} $ which will be of use in the sequel.
\blmma ( Section 3, \cite{majid_2} ) \label{30thsep20}
Let $\F$ and $\G$ be objects in the category $ {}_{\A} {\bf \mathcal{M}} {}_{\A}. $  The following statements hold:

{\rm (i)} We have a functor
$ {\rm bar}: {}_{\A} {\bf \mathcal{M}} {}_{\A} \rightarrow {}_{\A} {\bf \mathcal{M}} {}_{\A}  $ which sends an object $\F$ to $\overline{\F}$ and a morphism $\phi: \F \rightarrow \G$ to $ \overline{\phi}: \overline{F} \rightarrow \overline{G} $ defined by $ \overline{\phi} ( \overline{f} ) = \overline{\phi ( f )} $ for all $f$ in $\F.$

{\rm (ii)}  Viewing $\A$ as an object of $ {}_{\A} {\bf \mathcal{M}} {}_{\A}, $ we have a morphism $ \star: \A \rightarrow \overline{\A}, ~ a \mapsto \overline{a^*}. $

{\rm (iii)}  We have an invertible morphism in the category $ {}_{\A} {\bf \mathcal{M}} {}_{\A} $ given by
$$ \Upsilon: \overline{\F \tensora \G} \rightarrow \overline{\G} \tensora \overline{\F}, ~ \Upsilon ( \overline{f \tensora g} ) = \overline{g} \tensora \overline{f}. $$
$\Upsilon$ defines a natural equivalence between the functors $ {\rm bar} \circ \otimes $ and $ \otimes \circ ( {\rm bar} \times {\rm bar} ) \circ {\rm flip} $ from $ {}_{\A} {\bf \mathcal{M}} {}_{\A}  \times {}_{\A} {\bf \mathcal{M}} {}_{\A} $ to $ {}_{\A} {\bf \mathcal{M}} {}_{\A}.$

{\rm (iv)}  Define $ b b_{\F}: \F \rightarrow \overline{\overline{\F}} $ by the formula $ b b_{\F} ( f ) = \overline{\overline{f}}.$ Then $b b$ defines a natural equivalence between the functors $ \id $ and $ {\rm bar} \circ {\rm bar}. $ 
\elmma
This makes $ {}_{\A} {\bf \mathcal{M}} {}_{\A} $ into a bar-category ( see Definition 3.1 of \cite{majid_2} ). For more details on bar-categories, we refer to Subsection 2.8 of \cite{beggsmajidbook}. 

An object $\F$ of the bar category $ {}_{\A} {\bf \mathcal{M}} {}_{\A} $ is called a star object if there exists a morphism $ \star: \F \rightarrow \overline{\F} $ such that $ \overline{\star} \star ( f ) = \overline{\overline{f}} $ for all $f$ in $\F.$

For the rest of the section,  we will assume that $ (\E, d ) $ is a differential calculus on a $\ast$-algebra $\A$ such that $\E$ is a star object and moreover $ \star d = \overline{d} \star $ as maps from $\A$ to $\overline{\E}.$ Examples of such differential calculi have been studied in \cite{majid_2} and \cite{chern}. We explain one such example which will be used in the sequel.

\bxmpl \label{30thsep202}
Suppose $(\A, \mathcal{H}, D )$ is a spectral triple ( \cite{connes} ) over a unital $\ast$-algebra $\A,$ i.e, there exists a Hilbert space $\clh$ and a one-one $\ast$-homomorphism $ \pi: \A \rightarrow \clb ( \clh ) $ such that $ [ D, \pi ( a ) ] $ extends to a bounded operator on $\clh$ for all $a$ in $\A.$ Then, as in subsection 2.1.2 of \cite{frolich}, we have a canonical differential calculus $(\E, d )$ over $\A$ where the bimodule of one -forms 
$$\E = \{ \sum_{i} [ D, \pi ( a_i ) ] \pi ( b_i ) : a_i, b_i \in \A \} $$
and the restriction of the differential $d$ on $ \Omega^0 ( \A ) = \A $ is given by
$$ d: \A \rightarrow \E, ~ d ( a ) = \sqrt{-1} [ D, \pi ( a ) ]. $$
Since $d$ is a derivation, $\E$ is closed under  the adjoint operation $\ast$ inherited from $ \clb ( \clh ). $

We claim that $\E$ is a star-object in the bar category ${}_{\A} {\bf \mathcal{M}} {}_{\A}$ such that $ \star d = \overline{d} \star: \A \rightarrow \overline{\E}.$ 

Indeed, as in Subsection 2.2 of \cite{chern}, we define 
$$ \star ( e ) = \overline{e^*}, $$
where $e^*$ is the adjoint of the operator $e$ in $ \clb ( \clh ). $ Then for all $e$ in $\E,$
$$ \overline{\star} \star ( e ) = \overline{\star} ( \overline{e^*} ) = \overline{\overline{( e^* )^*}} = \overline{\overline{e}} $$
as the adjoint operation is an involution. Finally, for all $a$ in $\A,$
$$ \star d ( a ) = \star ( \sqrt{-1} [ D, \pi ( a ) ] ) = \overline{( \sqrt{-1} [ D, \pi ( a ) ]   )^*} = \overline{\sqrt{-1} [ D, \pi ( a^* ) ]} = \overline{d ( a^* )} = \overline{d} ( \overline{a^*} ) = \overline{d} \star ( a ).$$   
this finishes the proof of our claim.
\exmpl

 We have the following proposition from \cite{majid_2}.
 \bppsn ( \cite{majid_2} ) \label{2ndoct20}
Suppose $ (\E, d ) $ is a differential calculus over a star-algebra $\A$ such that the bimodule of one-forms $\E$ is a star-object in $ {}_{\A} {\bf \mathcal{M}} {}_{\A} $ and $ \star d = \overline{d} \star: \A \rightarrow \overline{\E}.$ If  $ ( \nabla, \sigma^\prime )$ is a right bimodule connection on an object  $\F$ in $ {}_{\A} {\bf \mathcal{M}} {}_{\A}, $ then we have a right bimodule connection $ ( \overline{\nabla}, ( \widehat{\sigma^\prime} )^{-1}   )$ on $\overline{\F},$ where 
$$ \overline{\nabla}: \overline{\F} \rightarrow  \overline{\F} \tensora \E, ~ \overline{\nabla} ( \overline{f} ) = ( \id \tensora \star^{-1} ) \Upsilon ( \overline{( \sigma^\prime)^{-1} \nabla ( f )} ) ~ {\rm and} $$
\begin{equation} \label{13thjan21} \widehat{\sigma^\prime}: \overline{\F} \tensora \E \rightarrow \E \tensora \overline{\F}: ~ \widehat{\sigma^\prime} = ( \star^{-1} \tensora \id ) \Upsilon \overline{\sigma^\prime} \Upsilon^{-1} ( \id \tensora \star ). \end{equation}
\eppsn
{\bf Proof:} Since the proof for left-connections is already available in \cite{majid_2}, we only provide a sketch of the proof. Indeed, following Section 3 of \cite{majid_2}, define
$$ \widehat{\nabla}: \overline{\F} \rightarrow \E \tensora \overline{\E}, ~ \widehat{\nabla} ( \overline{f} ) = ( \star^{-1} \tensora \id ) \Upsilon ( \overline{\nabla ( f )} ). $$
Then $ ( \widehat{\nabla}, \widehat{\sigma^\prime} ) $ can be checked to be a left bimodule connection on $\overline{\F}.$ Now, note that $ \overline{\nabla} = ( \widehat{\sigma^\prime} )^{-1} \widehat{\nabla}. $ Therefore, by the first assertion of Proposition \ref{27thsep20}, $ ( \overline{\nabla}, ( \widehat{\sigma^\prime} )^{-1} ) $ is a right bimodule connection on $ \overline{\F}. $
\qed

With a differential calculus $ (\E, d ) $ as in Proposition \ref{2ndoct20},  consider the category $ \C $ with objects $ ( \F, \nabla, \sigma^\prime ) $ where $\F$ is an $\A$-bimodule,  $ \sigma^\prime: \F \tensora \E \rightarrow \E \tensora \F $ is an $\A$-bilinear map and $ ( \nabla, \sigma^\prime ) $ is a right bimodule connection on $\F.$ A morphism between two objects $ ( \F, \nabla_1, \sigma^\prime )  $ and $ ( \mathcal{G}, \nabla_2, \tau   ) $ are $\A$-bimodule maps $ \theta: \F \rightarrow \mathcal{G} $ such that
$$ ( \theta \tensora \id ) \nabla_1 = \nabla_2 \circ \theta.  $$
Consider the map $ d: \A \rightarrow \A \tensora \E \cong \E, ~ d ( a ) = da. $ Then $ (\A, d, \id ) $ is an object of $\C.$ 
Using the second assertion of Proposition \ref{27thsep20}, it can be verified that $\C$ is a monoidal category with $ (\A, d, \id ) $ the identity object. Moreover, if $ (\F, \nabla, \sigma^\prime ) $ is an object of $\C$ and $ \widehat{\sigma^\prime}: \overline{\F} \tensora \E \rightarrow \E \tensora \overline{\F} $ as in \eqref{13thjan21}, then $ ( \overline{\E}, \overline{\nabla}, ( \widehat{\sigma^\prime} )^{-1} ) $ is another object of $\C.$ Moreover, the morphisms $\overline{\phi}, \Upsilon $ and $bb_{\F}$ are all morphisms in the category $\C.$

Beggs and Majid made the following definition.
\bdfn ( \cite{majid_2} )
Suppose $ ( \E, d ) $ is a differential calculus such that $\E$ is a star-object such that $ \star d = \overline{d} \star $ as maps from $\A$ to $\overline{\E}.$  A right bimodule connection  $ ( \nabla, \sigma^\prime ) $ on an $\A$-bimodule $\F$ is said to be star-compatible if
$$ ( \star \tensora \id ) \sigma^\prime ( e \tensora f ) = ( \widehat{\sigma^\prime} )^{-1} ( \id \tensora \star ) ( e \tensora f ) ~ {\rm for} ~ {\rm all} ~ e \in \E ~ {\rm and} ~ {\rm for} ~ {\rm all} ~ f \in \F. $$
Here, $ \widehat{\sigma^\prime} $ is as in \eqref{13thjan21}. 
\edfn 

\bthm 
Suppose $ ( \E, d ) $ is a tame differential calculus such that $\E$ is a star-object such that $ \star d = \overline{d} \star $ as maps from $\A$ to $\overline{\E}.$  Moreover, let $\sigma: \E \tensora \E \rightarrow \E \tensora \E $ denote the $\A$-bilinear map  as in Definition \ref{tame} Then the following statements hold:
\begin{enumerate}
 \item If $( \nabla, \sigma ) $ is a right bimodule connection on $\E,$ then $\nabla$ is star-compatible.

\item In particular, if $g$ is a {\bf bilinear} pseudo-Riemannian metric on $(\E, d ),$ then the Levi-Civita connection $\nabla$ for the triplet $ (\E, d, g ) $ is such that $(\nabla, \sigma ) $ is star-compatible.
\end{enumerate}  
\ethm 
{\bf Proof:} By part ii. of Lemma \ref{centeredremark}, it is enough to prove that for all $\omega, \eta \in \Ecenter $ and $a \in \A,$
\begin{equation} \label{28thsep20} ( \star \tensora \id ) \sigma ( \omega \tensora \eta a ) = ( \widehat{\sigma} )^{-1} ( \id \tensora \star ) ( \omega \tensora \eta a ). \end{equation}
Let $ \star ( \omega ) = \overline{e} $ and $ \star ( \eta ) = \overline{f}.$ As $\omega \in \Ecenter,$ $ b \omega = \omega b  $ for all $b$ in $\A$ and therefore, for all $b$ in $\A,$
$$ \overline{e b^*} = b \overline{e} = b \star ( \omega ) = \star ( b \omega ) = \star ( \omega b ) = \star ( \omega ) b = \overline{e} b  = \overline{b^* e }.$$
Hence, $ \overline{\overline{ e b^*}} = \overline{\overline{b^* e}}.$ Since ( by (iv) of Lemma \ref{30thsep20} ) we have a natural equivalence $bb_{\E},$ defined by $ bb_{\E} ( e ) = \overline{\overline{e}} $ from the functor $\id$ to the functor $ {\rm bar} \circ {\rm bar}, $ we deduce that for all $b$ in $\A,$ $ e b^* = b^* e, $ i.e, $e  $ belongs to $\Ecenter.$ We compute
\begin{eqnarray*} 
( \widehat{\sigma} )^{-1} ( \id \tensora \star ) ( \omega \tensora \eta a ) &=& ( \id \tensora \star^{-1} ) \Upsilon \overline{\sigma^{-1}} \Upsilon^{-1} ( \overline{e} \tensora \overline{f} a ) ~ {\rm (} ~ {\rm by} ~ \eqref{13thjan21} ~ {\rm )}\\
&=& ( \id \tensora \star^{-1} ) \Upsilon \overline{\sigma} ( \overline{f \tensora e}   a ) ~ {\rm (} ~ {\rm as} ~ \sigma^2 = \id ~ {\rm )}\\
&=& ( \id \tensora \star^{-1} ) \Upsilon \overline{\sigma ( f \tensora e )} a\\
&=& ( \id \tensora \star^{-1} ) \Upsilon ( \overline{e \tensora f} ) a ~ {\rm (} ~ {\rm by} ~ \eqref{10thjuly20182} ~ {\rm and} ~ {\rm since} ~ e  \in \Ecenter ~ {\rm )}\\
&=& ( \id \tensora \star^{-1} ) ( \overline{f} \tensora \overline{e} a )\\
&=& \star ( \eta ) \tensora \omega a\\
&=& ( \star \tensora \id ) \sigma ( \omega \tensora \eta a )
\end{eqnarray*}
by another application of \eqref{10thjuly20182}. This proves that $ ( \nabla, \sigma ) $ is star-compatible.

Now we prove 2. In Theorem 7.3 of \cite{article3}, it has been proven that if $g$ is bilinear, then  the Levi-Civita connection $ \nabla $ of $(\E, d, g )$ is such that $ (\nabla, \sigma ) $ is a right bimodule connection. Thus, the second assertion follows from the first one.
\qed

Stronger notions of compatibility with the star-structure have been studied in the literature. For the details, we  refer to Subsection 3.3 of \cite{majid_2} ( called $\ast$-preserving by the authors ) and Definition 3.2 of \cite{sitarz}. In this article, we will focus on the former one.

\bdfn ( eq. 7, \cite{majid_2} ) \label{2ndoct202}
Suppose $(\E, d )$ is a differential calculus over a $\ast$-algebra $\A$ as in Proposition \ref{2ndoct20}. A right bimodule connection $ (\nabla, \sigma^\prime ) $ on $\E$ is said to be star preserving if
\begin{equation} \label{1stoct203} \overline{\nabla} \star = ( \star \tensora \id ) \nabla \end{equation}
as maps from $\E$ to $\overline{\E} \tensora \E.$
\edfn	

If $\E$ is the bimodule of one-forms of the canonical differential calculus constructed out of a spectral triple, then from Example \ref{30thsep202}, we know that the differential calculus $(\E, d )$ satisfies the hypothesis of Proposition \ref{2ndoct20}. So Definition \ref{2ndoct202} makes sense for any bimodule connection on $\E.$ We will use this fact implicitly in the next two results.

\bthm \label{30thsep203}
Suppose $(\A, \clh, D )$ is a spectral triple over a unital $\ast$-algebra $\A$ such that the associated differential calculus $(\E, d )$ as in Example \ref{30thsep202} is tame. Let us continue to denote by $\sigma$ the canonical $\A$-bilinear map as in Definition \ref{tame} and let $(\nabla, \sigma )$ is a bimodule connection on the bimodule of one-forms $\E.$ Moreover, assume that
\begin{enumerate}
\item $\E$ is a free right $\A$-module generated by self-adjoint elements $e_i$ belonging to $\Ecenter.$

\item There exist self-adjoint elements $\Gamma^i_{j,k}$ such that $ \nabla ( e_i ) = \sum_{j,k} e_j \tensora e_k \Gamma^i_{j,k}.$
\end{enumerate}
Then the bimodule connection $(\nabla, \sigma )$ is star-preserving. 
\ethm
{\bf Proof:} As mentioned before, $\E$ is equipped with the adjoint operation inherited from $\clb ( \clh ) $ and $ \star: \E \rightarrow \overline{\E}  $ is given by $ \star ( e ) = \overline{e^*}. $ Moreover, since $\E$ is a free right $\A$-module generated by elements $ e_1, e_2, \cdots e_n  $ in $\Ecenter$, any element in $\E \tensora \E$ can be written as a unique linear combination $ \sum_{i,j} e_i \tensora e_j a_{ij} $ for some elements $a_{ij}$ in $\A.$ Therefore, that $\nabla ( e_i ) $ is of the form $  \sum_{j,k} e_j \tensora e_k \Gamma^i_{j,k}.$  So in (ii) of the statement, we are only assuming that $\Gamma^i_{j,k}$ is self-adjoint.

We begin by observing that it is enough to verify that for all $i = 1,2, \cdots n,$
\begin{equation} \label{1stoct20} \overline{\nabla} (  \star ( e_i ) ) = ( \star \tensora \id ) \nabla ( e_i ). \end{equation}
Indeed, for all $a$ in $\A,$  \eqref{1stoct20} implies that
\begin{eqnarray*}
\overline{\nabla} ( \star ( e_i ) a ) &=& \overline{\nabla} ( \overline{e^*_i} a )\\
&=& \overline{\nabla} ( \overline{e_i} a ) = \overline{\nabla} ( \overline{e_i}  ) a + \overline{e^*_i} \tensora da ~ {\rm (} ~ {\rm as} ~  e^*_i = e_i ~ {\rm )}\\
&=& \overline{\nabla} ( \star ( e_i ) ) a + \overline{e^*_i} \tensora da \\
&=& ( \star \tensora \id ) \nabla ( e_i  ) a + \overline{e^*_i} \tensora da\\
&=& ( \star \tensora \id ) ( \nabla ( e_i  ) a + e_i \tensora da )\\
&=& ( \star \tensora \id ) \nabla ( e_i a ).
\end{eqnarray*}
Since $ \{  e_i : i = 1,2, \cdots n \} $ is right $\A$-total in $\E,$ this finishes the proof of the theorem modulo the proof of \eqref{1stoct20}.

Now, in order to prove \eqref{1stoct20}, we observe that for all $i,$
\begin{equation} \label{1stoct202} ( {\Gamma^i_{j,k}} )^* \overline{e_i} = \overline{e_i \Gamma^i_{j,k}} = \overline{\Gamma^i_{j,k} e_i} = \overline{e_i} ( {\Gamma^i_{j,k}} )^* = \overline{e_i} \Gamma^i_{j,k} \end{equation}
as $ \Gamma^i_{j,k} $ is self-adjoint. We compute
\allowdisplaybreaks{
\begin{eqnarray*}
\overline{\nabla} ( \star ( e_i ) ) &=& \overline{\nabla} ( \overline{e_i} ) ~ {\rm (} ~ {\rm as} ~ e^*_i = e_i ~ {\rm )}\\
                                    &=& ( \id \tensora \star^{-1} ) \Upsilon ( \overline{\sigma \nabla ( e_i ) } ) ~ {\rm (} ~ {\rm as} ~ \sigma^2 = \id ~ {\rm )}\\
																		&=& ( \id \tensora \star^{-1} ) \Upsilon ( \overline{\sum_{j,k} e_k \tensora e_j \Gamma^i_{jk} } ) ~ {\rm (} ~ {\rm as} ~ e_j \in \Ecenter ~ {\rm and} ~ {\rm we} ~ {\rm have} ~ {\rm applied} ~ \eqref{10thjuly20182} ~ {\rm )}\\
																		&=& ( \id \tensora \star^{-1} ) ( \sum_{j,k} ( {\Gamma^i_{j,k}} )^* \overline{e_j} \tensora \overline{e_k} )\\
																		&=& ( \id \tensora \star^{-1} ) ( \sum_{j,k} \overline{e_j} {\Gamma^i_{j,k}} \tensora \overline{e_k} ) ~ {\rm (} ~ {\rm by} ~ \eqref{1stoct202} ~ {\rm )}\\	
																		&=& ( \id \tensora \star^{-1} ) ( \sum_{j,k} \overline{e_j} \tensora \Gamma^i_{j,k} \overline{e_k}  )\\
																		&=& ( \id \tensora \star^{-1} ) ( \sum_{j,k} \overline{e_j} \tensora  \overline{e_k \Gamma^i_{j,k}}  ) ~ {\rm (} ~ {\rm as} ~ \Gamma^i_{jk} ~ {\rm is} ~ {\rm self\text{-}adjoint} ~ {\rm)}\\
																		&=& ( \id \tensora \star^{-1} ) ( \sum_{j,k} \overline{e_j} \tensora  \overline{\Gamma^i_{j,k} e_k}  )\\
																		&=& \sum_{j,k} \overline{e_j} \tensora e_k \Gamma^i_{j,k}\\
																		&=& ( \star \tensora \id ) \nabla ( e_i )
\end{eqnarray*}
}
as $e^*_j = e_j$ and $ ( e_k \Gamma^i_{j,k} )^* = \Gamma^i_{j,k} e_k. $ This proves \eqref{1stoct20} and hence the theorem.
\qed

As an immediate corollary, we have:

\bcrlre
Consider the spectral triples on the noncommutative torus, quantum Heisenberg manifold and the Cuntz algebra on three generators as discussed in Subsection 6.1 of \cite{article7}, Section 6 of \cite{article1} and Section 3 of \cite{soumalya} respectively. Let $g_0$ be the unique pseudo-Riemannian bilinear metric on the bimodule of one-forms $\E$ as defined in Remark \ref{18thjan21}. Then the Levi-Civita connection for the triplet $(\E, d, g_0 )$ is star-preserving.
\ecrlre
{\bf Proof:} We need to verify that each of the differential calculi $(\E, d )$ and the Levi-Civita connection $\nabla$ satisfies the hypotheses of Theorem \ref{30thsep203}. The fact that these differential calculi are tame and that the bimodule one-forms are freely generated by self-adjoint central elements follow  from Proposition 6.10 of \cite{article7}, Proposition 6.3 ( and the definition of $e_i$ ) of \cite{article1} and Proposition 3.1 of \cite{soumalya}. 

In particular, by Remark \ref{18thjan21}, we indeed have a {\bf bilinear} pseudo-Riemannian metric $g_0$ on $\E$ such that $g_0 ( e_i \tensora e_j ) = \delta_{ij}. $ Since $g_0$ is a {\bf bilinear} pseudo-Riemannian metric, the Levi-Civita connection $\nabla$ for $(\E, d, g_0 )$ exists uniquely by Theorem \ref{21stfeb2020}. Moreover, $(\nabla, \sigma )$ is a right bimodule connection by virtue of Theorem 7.3 of \cite{article3}. Hence, the equation \eqref{1stoct203} makes sense.

Finally, we need to verify that the elements $\Gamma^i_{jk}$ appearing in (ii) of Theorem \ref{30thsep203} are self-adjoint. In all the above examples, $\Gamma^i_{jk}$ are real numbers as can be seen from the proof of Theorem 6.12 of \cite{article7} ( by putting $k = 1$ ), Theorem 6.17 of \cite{article7} and Theorem 4.4 of \cite{soumalya}. This completes the proof of the corollary.

\qed

		\section{Existence and uniqueness of Levi-Civita connection for strongly $\sigma$-compatible metrics} \label{exun}

In this section, we extend Theorem \ref{21stfeb2020} to the case of pseudo-Riemannian metrics which are strongly $\sigma$-compatible. By Definition \ref{metricdefn}, any pseudo-Riemannian metric $g$ on a tame differential calculus satisfies the compatibility equation $ g \circ \sigma = g. $ A strongly $\sigma$-compatible pseudo-Riemannian metric is a pseudo-Riemannian metric which satisfies an additional compatibility relation with $\sigma.$ 

Let $g$ be a pseudo-Riemannian metric on a tame differential calculus $ ( \E, d ) $ and $g^{(2)} : ( \E \tensora \E ) \tensora ( \E \tensora \E ) \rightarrow \A $ be defined by the following equation: 
$ g^{(2)}: ( \E \tensora \E ) \tensora ( \E \tensora \E ) \rightarrow \A $ is defined by 
\begin{equation} \label{3rdjune20} g^{(2)} ( ( e_1 \tensora f_1 ) \tensora ( e_2 \tensora f_2 )   ) = g ( e_1 g ( f_1 \tensora e_2 ) \tensora f_2 ). \end{equation}
\bdfn \label{stronglycompatible}
A pseudo-Riemannian metric $g$ on a tame differential calculus $( \E, d )$ is said to be strongly $\sigma$-compatible if for all $e_1, f_1, e_2, f_2  $ in $\E,$ the following equation holds:
\begin{equation} \label{3rdjune203} g^{(2)} ( \sigma ( e_1 \tensora f_1 ) \tensora ( e_2 \tensora f_2 )  ) = g^{(2)} (  ( e_1 \tensora f_1   ) \tensora \sigma ( e_2 \tensora f_2  )  ).  \end{equation}
\edfn

Now we prove a necessary and sufficient condition for a pseudo-Riemannian metric on a tame differential calculus to be strongly $\sigma$-compatible. 

\bppsn \label{3rdjune202}
A pseudo-Riemannian metric $g$ on a tame differential calculus $ ( \E, d ) $ is strongly $\sigma$-compatible if and only if any two elements of the set  $\{ g ( \omega \tensora \eta ): \omega, \eta \in \Ecenter   \} $ commute.

We have the following examples of strongly $\sigma$-compatible pseudo-Riemannian metrics:

\begin{enumerate}

\item Consider the canonical differential calculus on a manifold $M.$ Then any pseudo-Riemannian metric on $M$ is strongly $\sigma$-compatible.

\item Any bilinear pseudo-Riemannian metric on a tame differential calculus is strongly $\sigma$-compatible.

\item Suppose $( \E, d )$ is a tame differential calculus and $g_0$  a pseudo-Riemannian {\bf bilinear} metric on $\E.$ If $k$ is an invertible element of $\A,$ then the conformally deformed pseudo-Riemannian metric  $g = k. g_0 $ is  strongly $\sigma$-compatible.

\item Suppose $( \E, d )$ is a tame differential calculus such that the bimodule $\E$ of one-forms is finitely generated and free as a right $\A$-module. Let $ \{ e_1, e_2, \cdots e_n \} $ be a basis of $\E$ such that $e_i$ belongs to $\Ecenter$ for all $i.$  If $k$ is an invertible element of $\A,$ then the map 
$$ g:  \E \tensora \E \rightarrow \A, g ( ( \sum_i e_i a_i  ) \tensora ( \sum_j e_j b_j  ) ) = \sum_i k^i a_i b_i  $$
is a strongly $\sigma$-compatible pseudo-Riemannian metric on $\E.$
\end{enumerate}
\eppsn
{\bf Proof:} We begin by claiming that the equation \eqref{3rdjune203} holds if and only if for all $\omega_1, \eta_1, \omega_2, \eta_2 \in \Ecenter, $
\begin{equation} \label{3rdjune204} g^{(2)} (  \sigma ( \omega_1 \tensora \eta_1 ) \tensora ( \omega_2 \tensora \eta_2    ) ) = g^{(2)} ( ( \omega_1 \tensora \eta_1 ) \tensora \sigma( \omega_2 \tensora \eta_2  )     ).\end{equation}
Indeed, recall that by part iii. of Lemma \ref{centeredremark}, any element of $\E \tensora \E$ can be written as a linear combination $\sum_i v_i \tensora w_i a_i$ where $v_i, w_i \in \Ecenter $ and $a_i$ in $\A.$  If \eqref{3rdjune204} holds, then for all $v_i, w_i, v^\prime_j, w^\prime_j \in \Ecenter $ and $a_i, b_j \in \A, $ we get
\begin{eqnarray*}
g^{(2)} ( \sigma ( \sum_i v_i \tensora w_i a_i   ) \tensora ( \sum_j v^\prime_j \tensora w^\prime_j b_j    ) ) 
&=& \sum_{i,j} g^{(2)} (  \sigma ( v_i \tensora w_i ) \tensora ( v^\prime_j \tensora w^\prime_j  )  ) a_i b_j\\
&& {\rm(} ~ {\rm as} ~ \sigma ~ {\rm and} ~ g^{(2)} ~ {\rm are} ~ {\rm right} ~ \A-{\rm linear} ~ {\rm )}\\
&=& \sum_{i,j} g^{(2)} (  ( v_i \tensora w_i  ) \tensora \sigma ( v^\prime_j \tensora w^\prime_j ) ) a_i b_j\\
&=& g^{(2)} (  ( \sum_i v_i \tensora w_i a_i ) \tensora \sigma ( \sum_j v^\prime_j \tensora w^\prime_j b_j  )  ) 
\end{eqnarray*}
as $\sigma$ is bilinear by part iv. of Lemma \ref{17thdec20192}. This proves our claim.

Now if $\omega_1, \eta_1, \omega_2, \eta_2 \in \Ecenter $ and $g$ is strongly $\sigma$-compatible, then
\begin{eqnarray*}
g ( \omega_1 \tensora \eta_2 ) g ( \eta_1 \tensora \omega_2 ) &=& g ( \omega_1 \tensora \eta_2 g ( \eta_1 \tensora \omega_2 ) ) = g ( \omega_1 g ( \eta_1 \tensora \omega_2 ) \tensora \eta_2 ) ~ {\rm(} ~ {\rm as} ~ \eta_2 \in \Ecenter ~ {\rm )}\\
&=& g^{(2)} ( ( \omega_1 \tensora \eta_1  ) \tensora (  \omega_2 \tensora \eta_2   ) ) = g^{(2)} ( \sigma ( \eta_1 \tensora \omega_1   ) \tensora ( \omega_2 \tensora \eta_2 ) )\\
&=& g^{(2)} ( ( \eta_1 \tensora \omega_1  ) \tensora \sigma ( \omega_2 \tensora \eta_2  )  ) = g^{(2)} ( ( \eta_1 \tensora \omega_1    ) \tensora ( \eta_2 \tensora  \omega_2   ) )\\
&=& g ( \eta_1 \tensora \omega_2 ) g ( \omega_1 \tensora \eta_2 )
\end{eqnarray*}
and we have used that $g$ is right $\A$-linear, the equation \ref{17thdec20191} and $\omega_2$ belongs to $\Ecenter.$

Conversely, if any two elements of $ \{ g ( \omega \tensora \eta ): \omega, \eta \in \Ecenter    \} $ commute, then the above computation reveals that \eqref{3rdjune204} holds and therefore, by the claim made in the beginning of this proof, we deduce that $g$ is strongly $\sigma$-compatible.

Now we come to the examples. If $g$ is a pseudo-Riemannian metric on a manifold $M,$ then $ g ( e \tensora f ) $ belongs to $ C^\infty ( M ) = \mathcal{Z} ( C^\infty ( M ) )  $ and so it is strongly $\sigma = {\rm flip}$-compatible.
If $g$ is a pseudo-Riemannian bilinear metric on a tame differential calculus, then from Lemma 4.17 of \cite{article1}, we already know that $g$ is strongly $\sigma$-compatible. Otherwise, one can directly check it by using the first assertion along with part ii. of Lemma \ref{lemma0}.  

The case of the conformally deformed pseudo-Riemannian metric  is an easy consequence of part ii of Lemma \ref{lemma0}. If $ \omega, \eta, \omega^\prime, \eta^\prime \in \Ecenter $ and $g = k. g_0$ as in the statement, then 
$$ g ( \omega \tensora \eta ) g ( \omega^\prime \tensora \eta^\prime )  = k g_0 ( \omega \tensora \eta ) k g_0 ( \omega^\prime \tensora \eta^\prime ) = k g_0 ( \omega^\prime \tensora \eta^\prime ) k g_0 ( \omega \tensora \eta ) = g ( \omega^\prime \tensora \eta^\prime ) g ( \omega \tensora \eta ). $$
Hence, by the first assertion, $g$ is strongly $\sigma$-compatible.

Finally, we come to the fourth example. So let $  ( \E, d ) $  and $g$ be as in the statement. From Proposition \ref{11thaugust202}, we already know that $g$ is a pseudo-Riemannian metric. So we are left to prove that $g$ is strongly $\sigma$-compatible. We will again use the first assertion of this proposition. Let $\omega_1, \eta_1, \omega_2, \eta_2 $ belong to $\Ecenter. $ Then by Example \ref{11thaugust20}, we know that 
$$ \omega_1 = \sum_i e_i a_i, \eta_1 = \sum_j e_j b_j, \omega_2 = \sum_p e_p c_p, \eta_2 = \sum_q e_q d_q $$
for some elements $a_i, b_j, c_p, d_q$ belonging to $\Acenter.$ Hence,
\begin{eqnarray*}
g ( \omega_1 \tensora \eta_1 ) g ( \omega_2 \tensora \eta_2 ) &=& ( \sum_i k^i a_i b_i   ) ( \sum_p k^p c_p d_p  ) = \sum_i k^{i + p} a_i b_i c_p d_p \\
&=& ( \sum_p k^p c_p d_p  ) ( \sum_i k^i a_i b_i  ) = g ( \omega_2 \tensora \eta_2 ) g ( \omega_1 \tensora \eta_1 )
\end{eqnarray*}
as $ a_i, b_j c_p, d_q $ belong to $\Acenter.$ This completes the proof of the proposition.
 \qed

Now we are in a position to state the main result of the section.

\bthm
	 \label{existenceuniqueness}
	 If $ ( \E, d ) $ is a tame differential calculus and $g$ is a pseudo-Riemannian metric on $\E$ which is strongly $\sigma$-compatible, then in the presence of a pseudo-Riemannian {\bf bilinear} metric $g_0$ on $\E,$ there exists a unique Levi-Civita connection for the triplet $ ( \E, d, g ).$
	\ethm
	As an immediate corollary, we have the following:
	\bcrlre
	Consider the differential calculi on the fuzzy spheres, quantum Heisenberg manifold, Rieffel deformations of classical compact manifolds equipped with a free isometric toral action and Cuntz algebras of Example \ref{exampleassumption}. If $g$ is a strongly $\sigma$-compatible pseudo-Riemannian metric on any of these calculi, then there exists a unique Levi-Civita connection corresponding to $g.$
	\ecrlre
	{\bf Proof:} As remarked in Example \ref{exampleassumption}, all these calculi are tame. So we are left to verify the existence of a pseudo-Riemannian bilinear metric on each of these calculi. For the pseudo-Riemannian bilinear metric on the fuzzy spheres considered in \cite{article1} and \cite{article3}, we refer to Proposition 5.3  of \cite{article1} and Theorem 8.6 of \cite{article3} respectively.. Moreover, Proposition 6.4 of \cite{article1} and equation (6) of \cite{soumalya}  yield pseudo-Riemannian bilinear metrics on quantum Heisenberg manifolds and a Cuntz algebra respectively.   Finally, for the Rieffel deformation of classical compact manifolds as above, Proposition 7.25 of \cite{article1} guarantees the existence of such a metric.
	\qed
	 
 As mentioned in the introduction, we were unable to extend the Koszul formula based proof of \cite{article3} to the case of a general pseudo-Riemannian metric. Our approach is the one taken in \cite{article1}. Indeed, by Theorem \ref{19thfeb2020}, it is enough to prove that the map
\begin{equation} \label{29thmarch20} \Phi_g: {\rm Hom}_\A ( \E, \E \otimes^{{\rm sym}}_\A \E ) \rightarrow {\rm Hom}_\A ( \E \tensorsym \E, \E ) \end{equation}  
defined by
\begin{equation} \label{29thmarch202} \Phi_g ( L ) ( X ) = (g \tensora {\rm id}) \sigma_{23} (L \tensora {\rm id})(1 + \sigma) ( X ) \end{equation}
for all $ X \in \E \tensorsym \E $ is an isomorphism of right $ \A$-modules. Here, the maps $\sigma, \Psym, \E \tensorsym \E $ are as in Definition \ref{defnpsymandsigma}, Definition \ref{quasitame} and Definition \ref{19thfeb20202} respectively. 

Thus, throughout this section,  our goal will be to prove that the map $ \Phi_g $ is an isomorphism.
A crucial ingredient of the proof is the following proposition which was already proved in \cite{article1}:
\bppsn \label{psym23iso}
If $ ( \E, d ) $ is a tame differential calculus and the map $ \Psym: \E \tensora \E \rightarrow \E \tensora \E $ is the idempotent with range $ \E \tensorsym \E \subseteq \E \tensora \E $ ( and kernel $\F $ ) as defined in Definition \ref{quasitame}, then the map
$$ ( \Psym )_{23}:= ( {\rm id}_{\E} \tensora \Psym ): ( \E \tensorsym \E ) \tensora \E \rightarrow \E \tensora ( \E \tensorsym \E )  $$
is an isomorphism of right $\A$-modules. 
\eppsn
{\bf Proof:} We note that since $\Psym$ is left $\A$-linear by Remark \ref{17thdec2019remark}, the map $ ( \Psym )_{23}:= ( {\rm id}_{\E} \tensora \Psym ) $ is well-defined.

 Although this statement of this proposition is contained in \cite{article1}, let us sketch the proof for the sake of completeness. Since $( \E, d )$ is tame, $\E$ is centered by Lemma \ref{17thdec20192}. Then by Theorem 6.10 of \cite{Skd}, there exists a unique $\A$-bimodule isomorphism $\sigma^{{\rm can}}: \E \tensora \E \rightarrow \E \tensora \E$ such that  $(\sigma^{{\rm can}})^2= {\rm id}$  and  for all $\omega,\eta \in \Ecenter,$
$$\sigma^{{\rm can}}(\omega \tensora \eta)=\eta \tensora \omega.$$
Let $X = \frac{1 + \sigma^{{\rm can}}}{2}.$ Then by Lemma 2.8 of \cite{article1}, the map $ X_{23} $ is a right $\A$-linear isomorphism from $ {\rm Ran} ( X ) \tensora \E   $ onto $ X \tensora {\rm Ran} ( X ). $ 
Therefore, our proposition will be proved once we prove that our $\sigma$ is equal to $\sigma^{{\rm can}}.$ But the proof of this fact is contained in the proof of Proposition 6.3 of \cite{article3}.\qed

The proof of the isomorphism of \eqref{29thmarch20} will be derived by defining certain maps and proving some isomorphisms. In order to motivate our proof, let us recall the strategy employed in \cite{article1} for proving the isomorphism of $\Phi_g$ when $g$ is {\bf bilinear}. 

For a pseudo-Riemannian metric $g,$ let $V_g: \E \rightarrow \E^* $ be the map as in Definition \ref{metricdefn} and $g^{(2)}$ be as in \eqref{3rdjune20}. Thus,
$$ V_g ( e ) (  f ) = g ( e \tensora f )  $$
for all $e, f $ in $\E.$

Let us define $ V_{g^{(2)}}: \E \tensora \E \rightarrow ( \E \tensora \E )^* $ by the formula
$$ V_{g^{(2)}} ( \theta ) ( \theta^\prime ) = g^{(2)} ( \theta \tensora \theta^\prime ) ~ {\rm for} ~ {\rm all} ~ \theta,~ \theta^\prime ~ {\rm in} ~ \E \tensora \E. $$
When the pseudo-Riemannian metric $g$ is {\bf bilinear}, the isomorphism of the map $ \Phi_g $ follows from the following commutative diagram:
 \begin{equation} \label{diagram1} 
 \begin{tikzcd} 
  \Hom_{\A} ( \E, \E \tensorsym \E ) \arrow[r, "\zeta^{-1}_{\E \tensorsym \E, \E}"] \arrow[d, "\Phi_g"] & [5 ex] ( \E \tensorsym \E ) \tensora \E^* \arrow[r, "  \id \tensora V^{-1}_g"] & [5 ex] ( \E \tensorsym \E ) \tensora \E \arrow[d, "( \Psym  )_{23}"]\\
	\Hom_{\A} ( \E \tensorsym \E, \E ) &   \E \tensora ( \E \tensorsym \E )^* \arrow[swap, l, "\zeta_{\E, \E \tensorsym \E}"] &  \E \tensora ( \E \tensorsym \E ) \arrow[swap, l, "\id \tensora V_{g^{(2)}}"] 
 \end{tikzcd}
\end{equation}
For the definitions of $\zeta_{\E \tensorsym \E, \E}$ and $ \zeta_{\E, \E \tensorsym \E},$ we refer to Definition \ref{xi}.
Here,  $g$ is {\bf bilinear} and so the maps $ V_g $ and $g^{(2)}$ are $\A$-bilinear and so left $\A$-linear in particular. It follows that the map $V_{g^{(2)}}$ is also left $\A$-linear. Hence the maps $ \id \tensora V^{-1}_g  $ and $  \id \tensora V_{g^{(2)}} $ are well-defined and checked to be isomorphisms. Now, $ \zeta_{\E, \E \tensorsym \E} $ is a right $\A$-module isomorphism from $ \E \tensora ( \E \tensorsym \E   )^* $ to $ \Hom_\A ( \E \tensorsym \E, \E ) $ and by Proposition \ref{psym23iso},  $ ( \Psym )_{23} $ is an isomorphism from $ ( \E \tensorsym \E ) \tensora \E $ to $ \E \tensora ( \E \tensorsym \E ).$ Since the diagram \eqref{diagram1} is commutative,  $ \Phi_g $ is a composition of right $\A$-linear isomorphisms and hence an isomorphism. 

When $g$ is a general pseudo-Riemannian metric ( i.e, $g$ is only right $\A$-linear  ), the above argument fails. For example, the maps $ V_g $ and $V_{g^{(2)}}$ are only right $\A$-linear and so the maps $  \id \tensora V^{-1}_g  $ and $ \id \tensora V_{g^{(2)}}  $ are not well-defined. This is why the more general case of a strongly $\sigma$-compatible pseudo-Riemannian metric needs a different proof.
Our proof of the isomorphism of $\Phi_g$ in the general case follows from the following commutative diagram:
 \begin{equation} \label{diagram5}
 \begin{tikzcd} 
  \Hom_{\A} ( \E, \E \tensorsym \E ) \arrow[r, "\xi"] \arrow[d, "\Phi_g"] & ( \E \tensorsym \E ) \tensora \E \arrow[r, " 2 ( \Psym  )_{23}"] & [5 ex] \E \tensora (  \E \tensorsym \E ) \arrow[d, "\xi^{-1}"]\\
	\Hom_{\A} ( \E \tensorsym \E, \E ) & {} & \Hom_{\A} ( \E, \E \tensora \E ) \arrow[swap, ll, "\Psi_g"]
 \end{tikzcd}
\end{equation}
where  $\xi$ and $ \Psi_g $ are as defined in Definition \ref{xi2} and Definition \ref{defnpsig} respectively. We are going to prove ( Proposition \ref{13thapril202} ) that the maps $ \xi: \Hom_{\A} ( \E, \E \tensora \E ) \rightarrow \E \tensora \E \tensora \E$ and $ \Psi_g: \Hom_{\A} ( \E, \E \tensora \E ) \rightarrow \Hom_{\A} ( \E \tensora \E, \E ) $ are right $\A$-linear isomorphisms. Moreover, $\xi$ restricts to a right $\A$-linear isomorphism from $ \Hom_{\A} ( \E, \E \tensorsym \E ) $ to $ ( \E \tensorsym \E  ) \tensora \E $ ( Corollary \ref{19thjuly20184} ) while Proposition \ref{psym23iso} asserts that the map $\Psym$ is a right $\A$-linear isomorphism from $ ( \E \tensorsym \E ) \tensora \E $ to $ \E \tensora ( \E \tensorsym \E  ).$ So the commutativity of the diagram \eqref{diagram5} immediately implies that $ \Phi_g $ is one-one. The fact that $\Phi_g$ is onto $ \Hom_{\A} ( \E \tensorsym \E, \E ) $ also follows fairly easily once we prove the results mentioned above. 

\subsection{Some preparatory isomorphisms}

In this subsection, we define and prove some isomorphisms which will be needed in the proof of Theorem \ref{existenceuniqueness}. These include the maps $\xi$ and $\Psi_g$ appearing in  diagram \eqref{diagram5}. 

Throughout this section, we will work under the hypotheses of Theorem \ref{existenceuniqueness}. Thus, $ ( \E, d ) $ will denote a tame differential calculus, $g$ a strongly $\sigma$-compatible pseudo-Riemannian ( right $\A$-linear ) metric on $\E$ and   $g_0$ will denote a {\bf bilinear} pseudo-Riemannian bilinear metric on $\E.$ 

\brmrk \label{11thapril20}
By part iv. of Lemma \ref{lemma0}, $ V_g: \E \rightarrow \E^* $ is left $\Acenter$-linear. Thus, the map 
$$  \id \otimes_{\Acenter} V^{-1}_g: \Ecenter \otimes_{\Acenter} \E^* \rightarrow \Ecenter \otimes_{\Acenter} \E $$
 is well-defined. We will use this fact throughout the sequel.
\ermrk

Throughout the proof, certain isomorphisms will play a vital role which we note in the next lemma whose proof is elementary. 
\blmma \label{repeated}
Suppose $ ( \E, d ) $ is a tame differential calculus over $\A$ and $\mathcal{G}$ is an $\A$-bimodule. If $ u^{\E}: \Ecenter \otimes_{\Acenter} \A \rightarrow \E $ denotes the isomorphism in part i. of Definition \ref{tame}, then
$$ u^{\E} \tensora \id_{\mathcal{G}}: \Ecenter \otimes_{\Acenter} \mathcal{G} \cong ( \Ecenter \otimes_{\Acenter}  \A ) \tensora \mathcal{G}  \rightarrow \E \tensora \mathcal{G} $$
is left $\Acenter$-linear and right $\A$-linear isomorphism. If elements $ \omega_i  $ belong to $\Ecenter,$ $a_i$ belong to $\A$ and $f $ in $\mathcal{G},$ then
\begin{equation} \label{repeated1} ( u^{\E} \tensora \id_{\mathcal{G}}  ) ( \sum_i \omega_i \otimes_{\Acenter} a_i f   ) = \sum_i \omega_i a_i \tensora f.   \end{equation} 
Thus, the inverse map
$$ ( u^{\E} )^{-1} \tensora \id_{\mathcal{G}}: \E \tensora \mathcal{G} \rightarrow \Ecenter \otimes_{\Acenter} \mathcal{G} $$
is also a left $\Acenter$-linear and  right $\A$-linear isomorphism. If $ \omega_i, a_i $ and $f$ are as in \eqref{repeated1}, then 
\begin{equation} \label{repeated2} ( ( u^{\E} )^{-1} \tensora \id_{\mathcal{G}} ) ( \sum_i \omega_i a_i \tensora f ) = \sum_i \omega_i \otimes_{\Acenter} a_i f. \end{equation} 
\elmma
{\bf Proof:} Here and elsewhere, we are going to identify $ \A \tensora \mathcal{G} $ with $ \mathcal{G} $ without mentioning.  The equation \eqref{repeated1} follows directly from the definition of $u^\E.$ Since the proof of left $\Acenter$-linearity and right $\A$-linearity of $  u^{\E} \tensora \id_{\mathcal{G}}  $ is obvious, we omit their proofs.

 For the proof of \eqref{repeated2}, we observe that for $ \omega_i  $ in $\Ecenter$ and $ a_i $ in $\A,$ 
$$ ( u^\E )^{-1} ( \sum_i \omega_i a_i  ) = \sum_i \omega_i \otimes_{\Acenter} a_i $$
and so for all $f$ in $\mathcal{G},$
$$ ( ( u^\E )^{-1} \tensora {\rm id}  ) (  \sum_i \omega_i a_i \tensora f  ) = \sum_i \omega_i \otimes_{\Acenter} a_i \tensora f = \sum_i \omega_i \otimes_{\Acenter} a_i f. $$ 
\qed

We will need to use Lemma \ref{repeated} repeatedly in the sequel.

Now we define the ingredients $\xi$ and $ \Psi_g $ of diagram \eqref{diagram5} one by one. Firstly, in order to define $ \xi: \Hom_{\A} ( \E, \E \tensora \E ) \rightarrow \E \tensora \E \tensora \E,  $ we will need to introduce another
 map $ \tau: \E \tensora \E^* \rightarrow \E \tensora \E $ by the following commutative diagram: 
\begin{equation} \label{diagram2} 
\begin{tikzcd}
 \E \tensora \E^* \arrow[r, "( u^\E )^{-1} \tensora \id_{\E^*}"] \arrow[dr, "\tau"] & [5 ex] \Ecenter \otimes_{\Acenter} \E^* \arrow[r, "\id \otimes_{\Acenter} V^{-1}_g"] & [5 ex] \Ecenter \otimes_{\Acenter} \E \arrow[dl, " u^\E \tensora \id_{\E}" ] \\
{} & \E \tensora \E   & {}
\end{tikzcd}
\end{equation}
Here, the maps  $ ( u^{\E} )^{-1} \tensora \id_{\E^*}: \E \tensora \E^* \rightarrow \Ecenter \otimes_{\Acenter} \E^* $ and $ u^\E \tensora \id_{\E}: \Ecenter \otimes_{\Acenter} \E \rightarrow \E \tensora \E $ are the isomorphisms of Lemma \ref{repeated}. Moreover, the map $  \id \otimes_{\Acenter} V^{-1}_g: \Ecenter \otimes_{\Acenter} \E^* \rightarrow \Ecenter \otimes_{\Acenter} \E $ is well-defined by virtue of Remark \ref{11thapril20}. 

Using \eqref{repeated1} and \eqref{repeated2}, it is easy to check that if $\omega$ and $\omega^\prime$ belong to $\Ecenter,$ then
\begin{equation} \label{eqfortau} \tau ( \omega \tensora V_{g} ( \omega^\prime ) ) = \omega \tensora  \omega^\prime. \end{equation} 
It is clear that the map $\tau$ is left $\Acenter$-linear and so $ {\rm id} \otimes_{\Acenter} \tau $ is well-defined.
Now we can define the map $\xi.$
\bdfn \label{xi2}
We define $ \xi: \Hom_{\A} ( \E, \E \tensora \E ) \rightarrow \E \tensora \E \tensora \E  $ by the following commutative diagram:
\begin{equation} \label{diagram4}
\begin{tikzcd}
 \Hom_{\A} ( \E, \E \tensora \E ) \arrow[r, " \zeta^{-1}_{\E \tensora \E, \E}"] \arrow[ddrr, "\xi"] & [5 ex] ( \E \tensora \E  ) \tensora \E^* \arrow[r, "{( u^\E )}^{-1} \tensora \id_{\E \tensora \E^*}" ]  & [9 ex] \Ecenter \otimes_{\Acenter} ( \E \tensora \E^* ) \arrow[d, "\id \otimes_{\Acenter} \tau"]\\
 {} & {} & \Ecenter \otimes_{\Acenter} ( \E \tensora \E ) \arrow[d, "u^\E \tensora \id_{\E \tensora \E}"]\\
  {} & {} & \E \tensora \E \tensora \E
\end{tikzcd}
\end{equation}
\edfn
Here, $ ( u^{\E} )^{-1} \tensora \id_{\E \tensora \E^*}: \E \tensora ( \E \tensora \E^*  ) \rightarrow \Ecenter \otimes_{\Acenter} ( \E \tensora \E^*  ), $  $ u^{\E} \tensora \id_{\E \tensora \E}: \Ecenter \otimes_{\Acenter} ( \E \tensora \E  ) \rightarrow \E \tensora ( \E \tensora \E  ) $ are the isomorphisms as in Lemma \ref{repeated} and $ \zeta_{\E \tensora \E, \E}: ( \E \tensora \E ) \tensora \E^* \rightarrow \Hom_{\A} ( \E, \E \tensora \E ) $ is the right $\A$-linear isomorphism in Definition \ref{xi}.

Let $\omega, \eta, \eta^\prime$ be elements in $\Ecenter.$ The action of $\xi$  on the element $ \zeta_{\E \tensora \E, \E} ( \omega \tensora \eta \tensora V_{g} ( \eta^\prime )  )  $ can be computed from diagram \eqref{diagram4}. More precisely, using \eqref{repeated1}, \eqref{repeated2} and \eqref{eqfortau}, we obtain
\begin{equation} \label{13thapril20} \xi ( \zeta_{\E \tensora \E, \E} ( \omega \tensora \eta \tensora V_{g} ( \eta^\prime )  ) ) = \omega \tensora \eta \tensora \eta^\prime. \end{equation}
		
The definition of $\Psi_g$ also needs some preparations. 		
		\bppsn \label{hatinverse}
		Let $ ( \E, d ) $ be a tame differential calculus and $g$ a pseudo-Riemannian metric. We define ~ $ \widehat{\cdot}: ~ {\rm Hom}_{\A} ( \E, \E ) \rightarrow {\rm Hom}_{\A} ( \E \tensora \E, \A ) $ by $ B \mapsto \widehat{B,} $ where $ \widehat{B} $ is defined by:
$$ \widehat{B} ( \omega \tensora \eta ) = g ( B ( \omega ) \tensora \eta ). $$
Then the map
	$$ \Gamma: \E \tensora \E^* \rightarrow ( \E \tensora \E )^*:= {\rm Hom}_\A ( \E \tensora \E, \A ) ~~ {\rm defined} ~ {\rm by} ~~ \Gamma ( e \tensora \phi ) = \widehat{\zeta_{\E, \E} ( e \tensora \phi )}$$
	 is a left $ \mathcal{Z} ( \A ) $-linear  and a right $ \A $-linear isomorphism. Moreover, for all $e, f_1, f_2$ in $\E$ and $\phi$ in $\E^*,$
	 the following equation holds:
	\begin{equation} \label{15thapril20} \Gamma ( e \tensora \phi ) ( f_1 \tensora f_2 ) = g ( e \phi ( f_1 ) \tensora f_2 ). \end{equation}
\eppsn
{\bf Proof:} We start by claiming that the map $\widehat{\cdot}: \Hom_{\A} ( \E, \E ) \rightarrow \Hom_{\A} ( \E \tensora \E, \A ) $ is a left $\Acenter$-linear right $\A$-linear isomorphism. For $e, f$ in $\E$ and $ B $ in $\Hom_{\A} ( \E, \E ),$ we have
$$ \widehat{B a} ( e \tensora f ) = g ( B a ( e ) \tensora f  ) = g ( B ( a e ) \tensora f ) = \widehat{B} a ( e \tensora f ) $$
so that the map $\widehat{\cdot}$ is right $\A$-linear. Now, for $ a \in \Acenter $ and $e, f$ as above,
\begin{eqnarray*}
\widehat{ a B } ( e \tensora f ) &=& g ( a B ( e ) \tensora f )\\
                                 &=& g ( B ( e ) a \tensora f ) ~ {\rm (} ~ {\rm by} ~ \eqref{18thdec20191} ~  {\rm )}\\
																 &=& g ( B ( e ) \tensora a f )\\
																 &=& g ( B ( e ) \tensora f a ) ~ {\rm (} ~ {\rm by} ~ \eqref{18thdec20191} ~ {\rm )}\\
																 &=& g ( B ( e ) \tensora f ) a ~ {\rm (} ~ {\rm as} ~ g ~ {\rm is} ~ {\rm right} ~ \A\text{-}{\rm linear} ~ {\rm )}\\
&=& a g ( B ( e ) \tensora f ) ~ {\rm (} ~ {\rm as} ~ a \in \Acenter ~ {\rm )}\\
&=& a \widehat{B} ( e \tensora f )
\end{eqnarray*}
proving that $\widehat{\cdot}$ is left $\Acenter$-linear.

In order to prove that $ \widehat{\cdot} $ is an isomorphism, we directly define its inverse. Indeed, for $e$ in $\E$ and $C \in \Hom_{\A} ( \E \tensora \E, \A ), $ we define the element $ C ( e \tensora \cdot ) \in \E^* $ by the equation
 $$ C ( e \tensora \cdot ) ( f ) = C ( e \tensora f ). $$
Then we define a right $\A$-linear map 
$$ \cup: \Hom_{\A} ( \E \tensora \E, \A ) \rightarrow \Hom_{\A} ( \E, \E ), ~ C \mapsto C^{\cup,} $$
where 
 $$ C^\cup ( e ) := V^{-1}_g ( C ( e \tensora \cdot )  ).  $$
Thus, 
\begin{equation} \label{12thapril20} V_g ( C^\cup ( e )  ) = C ( e \tensora \cdot ). \end{equation}
Then for $  C \in \Hom_{\A} ( \E \tensora \E, \A ), $
$$ \widehat{ C^\cup } ( e \tensora f ) = g ( C^\cup ( e ) \tensora f ) = V_g ( C^\cup ( e ) ) ( f ) = C ( e \tensora f ) $$
by \eqref{12thapril20}.  

Similarly, if $ B $ belongs to $\Hom_{\A} ( \E, \E ), $ then by \eqref{12thapril20}, 
$$g ( ( \widehat{B} )^{\cup} ( e ) \tensora f  ) = V_g (  ( \widehat{B} )^{\cup} ( e )  ) ( f ) = \widehat{B} ( e \tensora f ) = g ( B ( e ) \tensora f )$$
for all $e, f$ in $\E.$
																			
By applying part iii. of Lemma \ref{lemma0}, this completes the proof that $ \cup $ is the inverse to $ \widehat{\cdot}. $

Let us recall ( Definition \ref{xi} ) that the map $ \zeta_{\E, \E} $ is a right $\A$-linear isomorphism from $ \E \tensora \E^* $ to $ \Hom_{\A} ( \E, \E ). $ It can be easily checked that it is also left $ \Acenter $-linear. We have proved that $ \widehat{\cdot} $ is left $\Acenter$-linear and right $\A$-linear isomorphism.   The map $\Gamma$ is the composition of $ \widehat{\cdot} $ with $ \zeta_{\E, \E} $ and so enjoys these properties.
Finally, the equation \eqref{15thapril20}  follows from the following computation:
\begin{eqnarray*}
\Gamma ( e \tensora \phi ) ( f_1 \tensora f_2 ) &=& \widehat{\zeta_{\E, \E} (  e \tensora \phi )} ( f_1 \tensora f_2 )\\
                                                &=& g ( \zeta_{\E, \E} ( e \tensora \phi ) ( f_1 ) \tensora f_2 )\\
																								&=& g ( e \phi ( f_1 ) \tensora f_2 ).
\end{eqnarray*}
\qed

Now we are in a position to define the map $ \Psi_g. $
\bdfn \label{defnpsig}
The map $\Psi_g: \Hom_{\A} ( \E, \E \tensora \E ) \rightarrow \Hom_{\A} ( \E \tensora \E, \E ) $ is defined via the following commutative diagram: 
\begin{equation} \label{diagram6}
 \begin{tikzcd}
  \Hom_{\A} ( \E, \E \tensora \E ) \arrow[r, "\zeta^{-1}_{\E \tensora \E, \E}"] \arrow[d, "\Psi_g"] & [5 ex] ( \E \tensora \E ) \tensora \E^* \arrow[r, "( u^\E )^{-1} \tensora \id_{\E \tensora \E^*}"] & [9 ex] \Ecenter \otimes_{\Acenter}  ( \E \tensora \E^* ) \arrow[d, "\id \otimes_{\Acenter} \Gamma"]\\
	\Hom_{\A} ( \E \tensora \E, \E ) &   \E \tensora ( \E \tensora \E )^* \arrow[swap, l, "\zeta_{\E, \E \tensora \E}"] & \Ecenter \otimes_{\Acenter} ( \E \tensora \E )^* \arrow[swap, l, "u^\E \tensora \id_{( \E \tensora \E )^*} "]
 \end{tikzcd}
\end{equation}
\edfn
Here, the maps $ \zeta_{\E \tensora \E, \E} $ and $ \zeta_{\E, \E \tensora \E} $ are as in Definition \ref{xi}.
 
Moreover, as $\Gamma$ is left $\Acenter$-linear by Proposition \ref{hatinverse}, the map $ {\rm id} \otimes_{\Acenter} \Gamma: \Ecenter \otimes_{\Acenter} ( \E \tensora \E^*  ) \rightarrow \Ecenter \otimes_{\Acenter} ( \E \tensora \E   )^* $ is well-defined.

Let  $\omega_1, \omega_2, \omega_3, \omega, \eta$ belong to $\Ecenter.$ Then the equation \eqref{repeated2} and the definition of $\zeta_{\E, \E \tensora \E}$ applied to diagram \eqref{diagram6} imply that 
\begin{eqnarray}
\Psi_g \zeta_{\E \tensora \E, \E} ( \omega_1 \tensora \omega_2 \tensora V_g ( \omega_3 ) ) ( \omega \tensora \eta ) &=& \omega_1  \Gamma ( \omega_2 \tensora V_g ( \omega_3 ) ) ( \omega \tensora \eta ) \nonumber \\
&=& \omega_1 g ( \omega_2 V_g ( \omega_3 ) ( \omega ) \tensora \eta ) ~ {\rm (} ~ {\rm by} ~ \eqref{15thapril20} ~ {\rm )} \nonumber \\
&=& \omega_1 g ( \omega_2 g ( \omega_3 \tensora \omega ) \tensora \eta ). \label{15thapril202}
\end{eqnarray}

  We end this subsection with the following result:
	
	\bppsn \label{13thapril202}
	The maps $\tau, \xi$ and $\Psi_g$ are all right $\A$-linear isomorphisms.
	\eppsn
	{\bf Proof:} For all the three maps, the isomorphism follows by showing that they are composition of right $\A$-linear isomorphisms.  From Lemma \ref{repeated}, we know that  $ ( u^{\E} )^{-1} \tensora {\rm id}_{\E^*} : \E \tensora \E^* \rightarrow \Ecenter \otimes_{\Acenter} \E^* $ and $ u^{\E} \tensora {\rm id}_{\E}: \Ecenter \otimes_{\Acenter} \E \rightarrow \E \tensora \E  $  are right $\A$-linear isomorphisms. By part ii. of Definition \ref{metricdefn}, $V^{-1}_g$ is a right $\A$-linear isomorphism from $\E^*$ to $\E.$ Thus, \eqref{diagram2} implies that $\tau$ is a composition of  right $\A$-linear isomorphisms.
	
From Definition \ref{xi}, we know that $ \zeta_{\E \tensora \E, \E} $ is a right $\A$-linear isomorphism from $ ( \E \tensora \E  ) \tensora \E^* $ to $ \Hom_{\A} ( \E, \E \tensora \E ) $ and so  	$\xi $ is an isomorphism of right $\A$-modules by \eqref{diagram4}.
	
	Finally, the isomorphism for $\Psi_g$ follows from the isomorphism property  of $ ( u^{\E} )^{-1} \tensora {\rm id}_{\E \tensora \E^*}, u^{\E} \tensora {\rm id}_{( \E \tensora \E  )^*} $ ( Lemma \ref{repeated} ), $ \zeta_{\E \tensora \E, \E} $ and $ \zeta_{\E, \E \tensora \E} $ ( Definition \ref{xi} ),  $\Gamma$ ( Proposition \ref{hatinverse} ) and  \eqref{diagram6}.
	\qed

\subsection{Proof of Theorem \ref{existenceuniqueness}}	

As explained at the end of the first subsection, the proof of Theorem \ref{existenceuniqueness} follows from the fact that  diagram \eqref{diagram5} is commutative. There are two main steps to prove that this diagram is commutative. These are Lemma \ref{basic} and Lemma \ref{phiandpsi}. We will also need to prove Corollary \ref{19thjuly20184} for proving that the map $\Phi_g$ is one-one and onto. Lemma \ref{23rdjuly20183} helps us to prove Lemma \ref{basic} and Lemma \ref{phiandpsi} while Lemma \ref{16thapril20} and  Lemma \ref{xilsigma} are used to prove Corollary \ref{19thjuly20184}.

\blmma \label{16thapril20}
The following statements hold:
\begin{enumerate}
\item If $g_0$ is the pseudo-Riemannian bilinear metric in Theorem \ref{existenceuniqueness}, then the sets $ S_1 = \{ \zeta_{\E \tensora \E, \E} ( \omega_1 \tensora \omega_2 \tensora V_{g_0} ( \omega_3 )  ): \omega_1, \omega_2, \omega_3 \in \Ecenter \} $ and $ S_2 =  \{ \zeta_{\E \tensora \E, \E} ( \omega_1 \tensora \omega_2 \tensora V_g ( \omega_3 )  ): \omega_1, \omega_2, \omega_3 \in \Ecenter \} $ are right $\A$-total in ${\rm Hom}_{\A} ( \E, \E \tensora \E ).$ 

\item If $ L $ belongs to $ \Hom_{\A} ( \E, \E \tensora \E ),$ then we define an element $ \sigma L $ belonging to $ \Hom_{\A} ( \E, \E \tensora \E ) $ as
 $$ ( \sigma L ) ( e ) = \sigma ( L ( e ) ). $$
$L$ belongs to $ \Hom_{\A} ( \E, \E \tensorsym \E ) $ if and only if $ \sigma L = L. $

\item If $ L $ belonging to $ \Hom_{\A} ( \E, \E \tensora \E ) $ is of the form $ \zeta_{\E \tensora \E, \E} ( \omega \tensora \eta \tensora V_{g} ( \omega^\prime ) ) $ for some $\omega, \eta, \omega^\prime $ in $\Ecenter,$ then
\begin{equation} \label{16thapril202} \sigma L = \zeta_{\E \tensora \E, \E} ( \eta \tensora \omega \tensora V_{g} ( \omega^\prime ) ). \end{equation} 
\end{enumerate}
\elmma
{\bf Proof:} The right $\A$-totality of $S_1$  has already been proved in part 3. of Lemma 4.16 of \cite{article1}. We will prove that $S_2$ is right $\A$-total in $ \Hom_{\A} ( \E, \E \tensora \E ) $ by using the right $\A$-totality of $S_1.$

Let $g_0$ be the pseudo-Riemannian {\bf bilinear} metric in Theorem \ref{existenceuniqueness}. Since $V_g$ is right $\A$-linear and $\Ecenter$ is right $\A$-total in $\E,$ the element
$ \omega_1 \tensora \omega_2 \tensora V_g ( V^{-1}_g V_{g_0} ( \omega_3 ) ) $ belongs to the right $\A$-linear span of $S_2.$ 

But
 $  \omega_1 \tensora \omega_2 \tensora V_g (  V^{-1}_g V_{g_0} ( \omega_3 ) ) = \omega_1 \tensora \omega_2 \tensora V_{g_0} ( \omega_3 ). $ 

Therefore,  
$ S_1  $ is contained in the right $\A$-linear span of $S_2.$ Since $ S_1 $ is right $\A$-total in $ \Hom_{\A} ( \E, \E \tensora \E ), $ the same is true about $S_2.$
This proves the first assertion.

 Now suppose $ L  $ belongs to $ \Hom_{\A} ( \E, \E \tensorsym \E ).$ Then for all $e$ in $\E,$ 
$$ ( \sigma L ) ( e ) = \sigma ( L ( e ) ) = ( 2 \Psym - 1  ) ( L ( e ) ) =  2 L ( e ) - L ( e ) =  L ( e ) $$
as $ L ( e ) $ belongs to the range of the idempotent $\Psym.$ The converse statement also follows easily.

Now we prove the third assertion. Using the definition of $\zeta_{\E \tensora \E, \E}$ ( Definition \ref{xi} ), it is easy to see that   for all $e$ in $\E,$ we have
\begin{equation} \label{2ndmay20} \zeta_{\E \tensora \E, \E} ( \omega \tensora \eta \tensora V_{g} ( \omega^\prime )  ) ( e ) = \omega \tensora \eta g ( \omega^\prime \tensora e ). \end{equation}
Hence, for all $e$ in $\E,$ we get
\begin{eqnarray*}
( \sigma L ) ( e ) &=& \sigma ( L ( e ) )\\
&=& \sigma ( \zeta_{\E \tensora \E, \E} ( \omega \tensora \eta \tensora V_{g} ( \omega^\prime )  ) ( e ) )\\
&=& \sigma ( \omega \tensora \eta g ( \omega^\prime \tensora e ) )\\
&=& \eta \tensora \omega g ( \omega^\prime \tensora e ) ~ {\rm (} ~ {\rm by} ~ \eqref{10thjuly20182} ~ {\rm )}\\
&=& \zeta_{\E \tensora \E, \E} ( \eta \tensora \omega \tensora V_{g} ( \omega^\prime )  ) ( e )
\end{eqnarray*}
and we have used the definition of $\zeta_{\E \tensora \E, \E}$ ( Definition \ref{xi} ) twice. Hence,
$$ \sigma L = \zeta_{\E \tensora \E, \E} ( \eta \tensora \omega \tensora V_{g} ( \omega^\prime )  ).$$
This completes the proof of the lemma.
\qed
		
\blmma \label{xilsigma}
Let $ L $ is an element of $ {\rm Hom}_\A ( \E, \E \tensora \E ), $ then
 \be \label{18thjuly20182} \xi (\sigma L) = \sigma_{12} \xi(L).\ee
 In particular,
\be \label{18thjuly20184} {\rm if} ~  L = \sigma L,~  {\rm then} ~  \xi(L )= \sigma_{12} \xi(L).\ee
\elmma
{\bf Proof:} We begin by claiming that \eqref{18thjuly20182} holds for $L$ of the form $ \zeta_{\E \tensora \E, \E} ( \omega \tensora \eta \tensora V_{g} ( \omega^\prime ) ) $ for all $\omega, \eta, \omega^\prime$ belonging to $\Ecenter.$ For such $L,$
\begin{eqnarray*}
\xi ( \sigma L ) &=& \xi ( \zeta_{\E \tensora \E, \E} (  \eta \tensora \omega \tensora V_{g} ( \omega^\prime )  ) ) ~ {\rm (} ~ {\rm by} ~ \eqref{16thapril202} ~ {\rm )}\\
&=& \eta \tensora \omega \tensora  \omega^\prime  ~ {\rm (} ~ {\rm by} ~ \eqref{13thapril20} ~ {\rm )}\\
&=& \sigma_{12} (  \omega \tensora \eta \tensora \omega^\prime   ) ~ {\rm (} ~ {\rm by} ~ \eqref{10thjuly20182} ~ {\rm )}\\
&=& \sigma_{12} \xi ( L )
\end{eqnarray*}
by \eqref{13thapril20}. This proves our claim.

Moreover, if  $ L^\prime \in \Hom_\A ( \E, \E \tensora \E ) $ satisfies \eqref{18thjuly20182} and $ a \in \A, $ then
 \begin{equation} \label{16thapril203} \xi ( \sigma (  L^\prime a  ) ) = \xi ( ( \sigma L^\prime ) a  ) = \xi (  \sigma L^\prime  ) a = ( \sigma_{12} \xi ( L^\prime )  ) a =  \sigma_{12} (  \xi ( L^\prime ) a   ) = \sigma_{12} \xi (  L^\prime a  ), \end{equation} 
where we have used that $ \sigma ( L a ) = (  \sigma L ) a $ and the right $\A$-linearity of the maps $ \sigma_{12} $ and $\xi.$  

Hence, if $ L^\prime \in \Hom_\A ( \E, \E \tensora \E ) $ satisfies \eqref{18thjuly20182}, so does the element $ L^\prime a. $ By the first assertion of Lemma \ref{16thapril20}, the set $ \{  \zeta_{\E \tensora \E, \E} ( \omega \tensora \eta \tensora V_{g} ( \omega^\prime ) ) : \omega, \eta, \omega^\prime \in \Ecenter \} $ is right $\A$-total in $ \Hom_\A ( \E, \E \tensora \E ). $ Therefore, \eqref{18thjuly20182} holds for all $ L  $ in $ \Hom_{\A} ( \E, \E \tensora \E ) $ by our claim and \eqref{16thapril203}.
\qed

\bcrlre \label{19thjuly20184}
 The map $\xi$ restricts to an isomorphism
 $$ \xi: \Hom_\A ( \E, \E \tensorsym \E ) \rightarrow ( \E \tensorsym \E ) \tensora \E. $$
\ecrlre
 {\bf Proof:} By Proposition \ref{13thapril202}, we already know that the map $ \xi $ is a right $\A$-linear isomorphism from $ \Hom_{\A} ( \E, \E \tensora \E ) $ to $ ( \E \tensora \E ) \tensora \E. $ Since $\E \tensorsym \E $ is a right $\A$-submodule of $\E \tensora \E,$ the map $\xi$ restricts to an isomorphism from $ \Hom_{\A} ( \E, \E \tensorsym \E ) $ onto its image. So we only need to check that $ \xi ( \Hom_{\A} ( \E, \E \tensorsym \E ) ) = ( \E \tensorsym \E  ) \tensora \E. $   

 Suppose $ L \in \Hom_\A ( \E, \E \tensorsym \E ).$ Then by part 2. of Lemma \ref{16thapril20}, $\sigma L = L. $ Hence, by \eqref{18thjuly20184}, we have $ \xi ( L ) = \sigma_{12} \xi ( L ). $ Thus, 
$$ ( \Psym )_{12} \xi ( L ) = \frac{\sigma_{12} + 1}{2} \xi ( L ) = \xi ( L ). $$
Therefore, $ \xi ( L ) \in ( \E \tensorsym \E ) \tensora \E. $

 Finally, if $ T \in ( \E \tensorsym \E ) \tensora \E \subseteq \E \tensora \E \tensora \E, $ Proposition \ref{13thapril202} implies the existence of a unique $ L^\prime \in \Hom_\A ( \E, \E \tensora \E ) $ such that $ \xi ( L^\prime ) = T. $ Since $  ( \Psym  )_{12} T = T $ ( so that $ \sigma_{12} T = T ), $ we have
$$ \xi ( L^\prime ) = T = \sigma_{12} T = \sigma_{12} \xi ( L^\prime ) = \xi ( \sigma L^\prime ) $$
by virtue of \eqref{18thjuly20182}. Since $ \xi $ is an isomorphism from $ \Hom_{\A} ( \E, \E \tensora \E ) $ to $ ( \E \tensora \E ) \tensora \E, $ we conclude that $ L^\prime = \sigma L^\prime, $ i.e, $ L^\prime \in \Hom_\A ( \E, \E \tensorsym \E ) $ again by part 2. of Lemma \ref{16thapril20}. This completes the proof of the corollary.
   \qed

	\blmma \label{23rdjuly20183}
	Let $ \omega_1, \omega_2, \omega_3 \in \Ecenter. $ Then we have:
	$$ \Psi_g \zeta_{\E \tensora \E, \E} (\omega_1 \tensora \omega_2 \tensora V_{g} ( \omega_3 ) ) \circ \sigma = \Psi_g \zeta_{\E \tensora \E, \E} (  \omega_1 \tensora  \omega_3  \tensora V_g ( \omega_2 ) ). $$
	\elmma
	{\bf Proof:} The proof follows by applying \eqref{15thapril202} and Proposition \ref{3rdjune202}. Let $\omega, \eta \in \mathcal{Z} ( \E ).$ Then  we obtain
					 \begin{eqnarray*}
					 &&\Psi_g \zeta_{\E \tensora \E, \E}  (\omega_1 \tensora \omega_2 \tensora V_{g} ( \omega_3 ) )(\sigma(\omega \tensora \eta))\\
					&=& \Psi_g \zeta_{\E \tensora \E, \E}  (\omega_1 \tensora \omega_2 \tensora V_{g} ( \omega_3 ) )( \eta \tensora \omega ) ~ {\rm (} ~ {\rm by} ~ \eqref{10thjuly20182} ~ {\rm )}\\
	        &=& \omega_1 g(\omega_2 g (\omega_3 \tensora \eta) \tensora \omega) ~ {\rm (} ~ {\rm by} ~ \eqref{15thapril202} ~ {\rm )}\\
                &=& \omega_1 g(\omega_2 \tensora \omega) g (\omega_3 \tensora \eta) ~ {\rm (} ~ {\rm as} ~ \omega \in \Ecenter ~ {\rm and} ~ g ~ {\rm is} ~ {\rm right} ~ \A\text{-}{\rm linear}   ~ {\rm )}\\
                &=& \omega_1 g (\omega_3 \tensora \eta)g(\omega_2 \tensora \omega) ~ {\rm (} ~ {\rm by} ~ {\rm Proposition} ~ \ref{3rdjune202} ~ {\rm )}\\
                &=& \omega_1 g( \omega_3 \tensora g(\omega_2 \tensora \omega) \eta) ~ {\rm (} ~ {\rm as} ~ \eta \in \Ecenter ~ {\rm )}\\
                &=& \Psi_g \zeta_{\E \tensora \E, \E}  (\omega_1 \tensora \omega_3 \tensora V_g(\omega_2))(\omega \tensora \eta)
     \end{eqnarray*}
			by \eqref{15thapril202}. Thus, for all $\omega, \eta$ in $\Ecenter,$ we have
$$ \Psi_g \zeta_{\E \tensora \E, \E} ( \omega_1 \tensora \omega_2 \tensora V_{g} ( \omega_3 ) ) ( \sigma ( \omega \tensora \eta ) ) = \Psi_g \zeta_{\E \tensora \E, \E} ( \omega_1 \tensora  \omega_3  \tensora V_{g} ( \omega_2 ) ) ( \omega \tensora \eta ). $$
		 Since  $ \{ \omega \tensora \eta: \omega, \eta \in \Ecenter   \} $ is right $\A$-total in $ \E \tensora \E $ by part ii. of Lemma \ref{centeredremark} and the maps  $\Psi_g \zeta_{\E \tensora \E, \E} ( \omega_1 \tensora \omega_2 \tensora V_g ( \omega_3 ) ), \Psi_g \zeta_{\E \tensora \E, \E} ( \omega_1 \tensora \omega_3 \tensora V_g ( \omega_2 ) ) $ are right $\A$-linear by Proposition \ref{13thapril202}, the proof is complete.
      \qed

From Proposition \ref{13thapril202}, we know that $\Psi_g$ is a right $\A$-linear isomorphism from $ \Hom_{\A} ( \E, \E \tensora \E ) $ to $ \Hom_{\A} ( \E \tensora \E, \E ). $ This allows us to make the next definition. 

\bdfn \label{13thapril204}
We define the map $ \Theta_g: \Hom_{\A} ( \E \tensora \E, \E ) \rightarrow \Hom_{\A} ( \E, \E \tensora \E ) $ to be the inverse of $\Psi_g.$
\edfn

	 \blmma
	 \label{basic}
	Let $L \in {\rm Hom}_{\A}(\E, \E \tensora \E)$ and $\Theta_g$ be as in Definition \ref{13thapril204}. Then
	 \be \label{19thjuly20182} \xi(  \Theta_g ( \Psi_g ( L ) \circ \sigma   ) ) = \sigma_{23} \xi(L). \ee
	 \elmma
	 {\bf Proof:} We begin by claiming that it is enough to prove \eqref{19thjuly20182} for $ L $ belongs to the set $S_2$ defined in Lemma \ref{16thapril20}. 
	Indeed, if $L \in {\rm Hom}_{\A}(\E, \E \tensora \E)$ satisfies \eqref{19thjuly20182} and $ l_a $ denotes the left multiplication by the element $a,$ then by using the right $\A$-linearity of $ \sigma, \xi, \Theta_g, \Psi_g $ and the left $\A$-linearity of $ \sigma $ ( the fourth assertion of  Lemma \ref{17thdec20192}  ),   we have:
	\begin{eqnarray*}
	  \sigma_{23} \xi ( L a ) &=& (  \sigma_{23} \xi ( L )   ) a = ( \xi (  \Theta_g (  \Psi_g ( L ) \circ \sigma   )  ) ) a\\
	 &=& \xi (  \Theta_g (  \Psi_g ( L ) \circ \sigma   ) a ) = \xi ( \Theta_g (  \Psi_g ( L ) \circ \sigma \circ l_a   )    )\\
	&=& \xi (  \Theta_g (  \Psi_g ( L ) \circ l_a \circ \sigma  )  ) = \xi (  \Theta_g (  \Psi_g (  La  ) \circ \sigma  )  ),
	\end{eqnarray*}
	i.e, $ La $ also satisfies \eqref{19thjuly20182}. 
	Since the set $S_2$ is right $\A$-total in $\Hom_{\A} ( \E, \E \tensora \E ) $ by the first assertion of Lemma \ref{16thapril20}, this proves our claim.

	Now we prove that the equation \eqref{19thjuly20182} indeed holds for $L$ belonging to the set $S_2.$
	
So let $ L = \zeta_{\E \tensora \E, \E} ( \omega_1 \tensora \omega_2 \tensora V_g ( \omega_3 ) ) $ for some  $ \omega_1, \omega_2, \omega_3 \in \Ecenter. $ Then by Lemma \ref{23rdjuly20183}, we have,
	\begin{eqnarray*}
	 &&\xi (  \Theta_g ( \Psi_g \zeta_{\E \tensora \E, \E} ( \omega_1 \tensora \omega_2 \tensora V_g (  \omega_3    )   ) \circ \sigma  )  )\\
	 &=& \xi (  \Theta_g  \Psi_g \zeta_{\E \tensora \E, \E} ( \omega_1 \tensora   \omega_3   \tensora V_g ( \omega_2 )  )   )\\
	 &=& \xi \zeta_{\E \tensora \E, \E} ( \omega_1 \tensora   \omega_3   \tensora V_g ( \omega_2 )  )\\
	 &=& \omega_1 \tensora  \omega_3  \tensora \omega_2 ~ {\rm (} ~ {\rm by} ~ \eqref{13thapril20}   ~ {\rm )}\\
	 &=& \sigma_{23} ( \omega_1 \tensora \omega_2 \tensora \omega_3  )\\
	 &=& \sigma_{23} \xi ( L ).  ~ (~ {\rm by} ~ \eqref{13thapril20}   ~)\\
	\end{eqnarray*}\qed

    \blmma \label{phiandpsi}
    For $L \in {\rm Hom}_{\A} (\E, \E \tensora^{\rm sym} \E)$ and $ X \in \E \tensorsym \E, $ we have 
	 \be \label{19thjuly2018}	\Phi_g(L) ( X ) = ( \Psi_g ( L) \circ (1+\sigma) ) ( X ). \ee
    \elmma
{\bf Proof:} By part 1 of Lemma \ref{16thapril20}, we know that the set $ S_2 =  \{ \zeta_{\E \tensora \E, \E} ( \omega_1 \tensora \omega_2 \tensora V_{g} ( \omega_3 )  ): \omega_1, \omega_2, \omega_3 \in \Ecenter \} $ is right $\A$-total in $ \Hom_{\A} ( \E, \E \tensora \E ). $ Moreover, $ \{ \omega \tensora \eta: \omega, \eta \in \Ecenter \} $ is right $\A$-total in $ \E \tensora \E $ by Lemma \ref{centeredremark}   and so $ \{ \Psym ( \omega \tensora \eta ) : \omega, \eta \in \Ecenter \} $ is right $\A$-total in $ \E \tensorsym \E = {\rm Ran} ( \Psym ). $ So, it suffices to check that for all $ L \in S_2, \omega, \eta \in \Ecenter  $ and $a, b \in \A,$ the following equation holds: 
	\begin{equation} \label{3rdmay20} \Phi_g ( L a ) ( \Psym ( \omega \tensora \eta b ) ) = \Psi_g ( L a ) ( 1 + \sigma ) \Psym ( \omega \tensora \eta b ). \end{equation}
	But this easily follows once we check that for all $L, \omega, \eta $ as above,
	\begin{equation} \label{3rdmay202}  \Phi_g ( L  ) ( \Psym ( \omega \tensora \eta  ) ) = \Psi_g ( L  ) ( 1 + \sigma ) \Psym ( \omega \tensora \eta  ).  \end{equation}
	Indeed, we use the right $\A$-linearity of the maps $ \Phi_g, \Psi_g, \Psym $ multiple times to compute
	\begin{eqnarray*}   
	\Phi_g ( L a ) ( \Psym ( \omega \tensora \eta b ) ) &=& \Phi_g ( L ) a \Psym ( \omega \tensora \eta ) b\\ 
	&=&  \Phi_g ( L ) \Psym ( a \omega \tensora \eta ) b ~ {\rm (}  ~ {\rm as} ~ \Psym = \frac{1 + \sigma}{2} ~ {\rm is} ~ {\rm left} ~ \A\text{-}{\rm linear} ~ {\rm by} ~ {\rm Lemma} ~ \ref{17thdec20192} ~ {\rm )}\\
	&=& \Phi_g ( L ) \Psym ( \omega \tensora \eta ) a b ~ {\rm (} ~ {\rm as} ~ \omega, \eta \in \Ecenter ~ {\rm )}\\
	&=& \Psi_g ( L ) ( 1 + \sigma ) \Psym ( \omega \tensora \eta ) a b ~ {\rm (} ~ {\rm by} ~ \eqref{3rdmay202} ~ {\rm )}\\
	&=&  \Psi_g ( L a ) ( 1 + \sigma ) \Psym ( \omega \tensora \eta b )  
\end{eqnarray*}
since $\Psym$ and $\sigma$ are left $\A$-linear. This proves \eqref{3rdmay20}. So we only need to prove \eqref{3rdmay202}. If $L = \zeta_{\E \tensora \E, \E} ( \omega_1 \tensora \omega_2 \tensora V_g ( \omega_3 ) ) $ for some $\omega_1, \omega_2, \omega_3 \in \Ecenter, $ then
\begin{eqnarray*}
\Phi_g ( L ) \Psym ( \omega \tensora \eta ) &=&  ( g \tensora {\rm id} ) \sigma_{23} ( L \tensora {\rm id} ) ( 1 + \sigma ) \frac{1 + \sigma}{2} ( \omega \tensora \eta ) ~ {\rm (} ~ {\rm by} ~ \eqref{29thmarch202} ~ {\rm )}\\
&=& ( g \tensora {\rm id} ) \sigma_{23} ( L \tensora {\rm id} ) ( \omega \tensora \eta + \eta \tensora \omega )\\
&& ~ {\rm(} ~ {\rm as} ~ \Psym = \frac{1 + \sigma}{2} ~ {\rm is} ~ {\rm an} ~ {\rm idempotent} ~ {\rm and} ~ {\rm we} ~ {\rm have} ~ {\rm used} ~ \eqref{10thjuly20182} ~ {\rm)}\\
&=& ( g \tensora {\rm id} ) \sigma_{23} ( L ( \omega ) \tensora \eta + L ( \eta ) \tensora \omega )\\
&=& ( g \tensora {\rm id} ) \sigma_{23} ( \omega_1 \tensora \omega_2 g ( \omega_3 \tensora \omega ) \tensora \eta + \omega_1 \tensora \omega_2 g ( \omega_3 \tensora \eta ) \tensora \omega )   ~ {\rm (} ~ {\rm by} ~ \eqref{2ndmay20} ~ {\rm)}\\  
&=& ( g \tensora {\rm id} ) ( \omega_1 \tensora \eta \tensora \omega_2 g ( \omega_3 \tensora \omega ) + \omega_1 \tensora \omega \tensora \omega_2 g ( \omega_3 \tensora \eta ) )  ~ {\rm (} ~ {\rm by} ~ \eqref{10thjuly20182} ~ {\rm )}\\
&=& g ( \omega_1 \tensora \eta ) \omega_2 g ( \omega_3 \tensora \omega ) + g ( \omega_1 \tensora \omega ) \omega_2 g ( \omega_3 \tensora \eta )\\
&=& \omega_2 g ( \omega_1 \tensora \eta ) g ( \omega_3 \tensora \omega ) + \omega_2 g ( \omega_1 \tensora \omega ) g ( \omega_3 \tensora \eta )
\end{eqnarray*}
as $\omega_2$ belongs to $\Ecenter.$

Now, the equation \eqref{16thapril202} implies that $ \sigma L = \zeta_{\E \tensora \E, \E} ( \omega_2 \tensora \omega_1 \tensora V_{g} ( \omega_3 ) ) $ and hence by \eqref{15thapril202}, we obtain
\begin{eqnarray*}
&& \omega_2 g ( \omega_1 \tensora \eta ) g ( \omega_3 \tensora \omega ) + \omega_2 g ( \omega_1 \tensora \omega ) g ( \omega_3 \tensora \eta )\\
&=& \Psi_g ( \sigma L ) ( \omega \tensora \eta + \eta \tensora \omega )\\
&=& \Psi_g ( L ) ( \omega \tensora \eta + \eta \tensora \omega ) ~ {\rm (} ~ {\rm as} ~ \sigma L = L ~ {\rm by} ~ {\rm part} ~ {\rm 2.} ~ {\rm of} ~ {\rm Lemma} ~ \ref{16thapril20} ~ {\rm)}\\
&=& \Psi_g ( L ) ( 1 + \sigma ) ( \omega \tensora \eta ) ~ {\rm (} ~ {\rm by} ~ \eqref{10thjuly20182} ~ {\rm )}\\
&=& \Psi_g ( L ) 2. \frac{( 1 + \sigma )}{2} \frac{( 1 + \sigma )}{2} ( \omega \tensora \eta ) ~ {\rm (} ~ {\rm as} ~ \frac{1 + \sigma}{2} = \Psym ~  {\rm is} ~ {\rm an} ~ {\rm idempotent} ~ {\rm )}\\
&=& \Psi_g ( L ) ( 1 + \sigma ) \Psym ( \omega \tensora \eta ).  
\end{eqnarray*}
This completes the proof of \eqref{3rdmay202} and hence the proposition.
\qed

	Now we have all the ingredients to prove that the diagram \eqref{diagram5} commutes.
	\bppsn \label{13thapril203}
	For all $L$ in $\Hom_{\A} ( \E, \E \tensorsym \E ) $ and for all $X $ in $ \E \tensorsym \E, $ the following equation holds:
	\begin{equation} \label{23rdjuly2018} \Phi_g ( L ) ( X ) = ( \Psi_g \xi^{-1} 2 ( \Psym)_{23} \xi ( L )  ) ( X ). \end{equation}
	\eppsn
	{\bf Proof:} The proof follows by a combination of Lemma \ref{phiandpsi} and Lemma \ref{basic}. Indeed, we get
	\begin{eqnarray*}
	\Phi_g ( L ) ( X ) &=& \Psi_g ( L ) ( 1 + \sigma ) ( X ) ~ {\rm(} ~ {\rm by} ~ {\rm Lemma} ~ \ref{phiandpsi} ~ {\rm)}\\
	&=& \Psi_g ( L ) ( X ) + \Psi_g ( L ) \circ \sigma ( X )\\
	&=& \Psi_g ( L ) ( X ) + \Psi_g \xi^{-1} \sigma_{23} \xi ( L ) ( X ) ~ {\rm(} ~ {\rm by} ~ {\rm Lemma} ~ \ref{basic} ~ {\rm)}\\
	&=& ( \Psi_g \xi^{-1} \xi ( L ) + \Psi_g \xi^{-1} \sigma_{23} \xi ( L ) ) ( X )\\
	&=& \Psi_g \xi^{-1} ( 1 + \sigma_{23} ) \xi ( L ) ( X )\\
	&=& \Psi_g \xi^{-1} 2 ( \Psym )_{23} \xi ( L ) ( X ).
	\end{eqnarray*}
\qed
	
	\vspace{4mm}

We are now in a position to prove the main result.

\vspace{4mm}
	 
{\bf Proof of Theorem \ref{existenceuniqueness}:} As explained at the end of the first subsection of this section, the fact that $\Phi_g$ is one-one immediately follows from \eqref{diagram5} ( now proved in Proposition \ref{13thapril203} ). So we are left to prove that $ \Phi_g $ is onto.

Let $ T $ be an element of $ \Hom_\A ( \E \tensorsym \E, \E ). $ Let us define $ \widetilde{T} \in \Hom_\A ( \E \tensora \E, \E ) $ by the formula
  $$ \widetilde{T} = T \circ \Psym. $$
	Let $\Theta_g: \Hom_{\A} ( \E \tensora \E, \E ) \rightarrow \Hom_{\A} ( \E, \E \tensora \E ) $ be the inverse of the map $ \Psi_g $ as defined in Definition \ref{13thapril204}. We claim that 
$$ L:= \frac{1}{2} \xi^{-1} (   \Psym   )^{-1}_{23} \xi \Theta_g ( \widetilde{T} ) $$
is well defined, belongs to $ \Hom_\A ( \E, \E \tensorsym \E ) $ ( i.e, $ {\rm Ran} ( L ) \subseteq \E \tensorsym \E $ ) and $ \Phi_g ( L ) = T. $

Since $ \Psi_g $ is an isomorphism from $ \Hom_{\A} ( \E, \E \tensora \E ) $ to $ \Hom_{\A} ( \E \tensora \E, \E ) $ by Proposition \ref{13thapril202}, there exists a unique $ W \in \Hom_\A ( \E, \E \tensora \E ) $ such that $ \Psi_g ( W ) = \widetilde{T}.  $ Therefore,
\begin{eqnarray*}
\xi ( \Theta_g (  \Psi_g ( W ) \circ \sigma  ) ) &=& 
\xi ( \Theta_g (  \widetilde{T} \circ \sigma  ) )\\
 &=& \xi ( \Theta_g (  \widetilde{T}  ) ) ~  ( ~   {\rm since} ~ \widetilde{T} = T \circ \Psym ~ {\rm and} ~ \sigma = 2 \Psym - 1   ~ )\\
&=& \xi ( \Theta_g \Psi_g ( W ) )\\
&=& \xi ( W )
\end{eqnarray*}
as $ \Theta_g $ is the invserse of $ \Psi_g. $
Hence, by \eqref{19thjuly20182}, we have $  \xi ( W ) = \sigma_{23} \xi ( W ) $ and thus, 
$$ \xi ( W ) = \frac{1 + \sigma_{23}}{2} \xi ( W ) = (  \Psym )_{23} \xi ( W ). $$
Consequently, $ \xi \Theta_g ( \widetilde{T} ) = \xi ( W ) $ belongs to $ {\rm Ran} ( (  \Psym  )_{23} ) $ and so by Proposition \ref{psym23iso}, $ ( \Psym )^{-1}_{23} ( \xi \Theta_g ( \widetilde{T} ) )  $ exists and belongs to $ (   \E \tensorsym \E ) \tensora \E.  $ Since $\xi$ is an isomorphism from $\Hom_{\E} ( \E, \E \tensorsym \E ) $ to $ ( \E \tensorsym \E ) \tensora \E  $ ( Corollary \ref{19thjuly20184} ),  $ L =  \frac{1}{2} \xi^{-1} (   \Psym   )^{-1}_{23} ( \xi \Theta_g ( \widetilde{T} ) ) $ belongs to $ \Hom_\A ( \E, \E \tensorsym \E ).$

   Finally, \eqref{23rdjuly2018} shows that for all $ X $ in $\E \tensorsym \E,$ 
	$$ \Phi_g ( L ) ( X ) = (  \Psi_g \xi^{-1} 2 ( \Psym  )_{23} \xi \frac{1}{2} \xi^{-1} ( \Psym )^{-1}_{23} \xi \Theta_g \widetilde{T}   ) ( X ) = ( \Psi_g \Theta_g \widetilde{T}  ) ( X ) = \widetilde{T} ( X ) = T \circ \Psym ( X ) = T ( X )  $$
	as $\Psym$ is an idempotent onto $ \E \tensorsym \E $ and $ X $ belongs to $\E \tensorsym \E.$  
	Hence, $ \Phi_g $ is onto.
	\qed
	
	\vspace{4mm}

{\bf Acknowledgement}  The authors would like to thank
 Ulrich Krahmer for several useful comments and discussion, in particular pointing out to them the reference \cite{heckenberger_etal}. D.G will like to thank D.S.T, Government of India for J.C. Bose National Fellowship. S.J will like to thank D.S.T, Government of India for the Inspire Fellowship.


\begin{thebibliography}{6666}

\bibitem{pseudo} J. Arnlind and M. Wilson: Riemannian curvature of the noncommutative 3-sphere, J. Noncomm. Geom. {\bf 11} (2017)  507--536.

\bibitem{cylinder} J. Arnlind and G. Landi:  Projections, modules and connections for the noncommutative cylinder, Adv. Theor. Math. Phys. {\bf 24} (2020) 527--562.

\bibitem{tiger} J. Arnlind and A.T. Norkvist: Noncommutative minimal embeddings and morphisms of pseudo-Riemannian calculi, J. Geom. Phys. {\bf 159} (2020).

\bibitem{landiqhom} J. Arnlind, K. Ilwale and G. Landi: On q-deformed Levi-Civita connections, arXiv: 2005.02603.

\bibitem{aschieri} P. Aschieri: Cartan structure equations and Levi-Civita connection in braided geometry, arXiv:2006.02761.


	
 \bibitem{majid_2} E.J. Beggs and S. Majid: $\ast$-Compatible connections in noncommutative Riemannian geometry, J. Geom. Phys. {\bf 61} (2011) 95--124.

\bibitem{chern} E. Beggs and S. Majid: Spectral triples from bimodule connections and Chern connections, Volume {\bf 11}, Issue 2, 2017, pp. 669--701.

\bibitem{beggsmajidbook} E.J. Beggs and  S. Majid: Quantum Riemannian geometry, Grundlehren der mathematischen Wissenschaften, Springer Verlag, 2019. 

	
\bibitem{article1} J. Bhowmick, D. Goswami and S. Mukhopadhyay: Levi-Civita connections for a class of spectral triples, Letters in Mathematical Physics, {\bf 110} ( 2020 ), 835--884. 
			
			
			
 \bibitem{article3} J. Bhowmick, D. Goswami and G. Landi: On the Koszul formula in noncommutative geometry, {Reviews in Mathematical Physics}, {\bf 32}, No 10, 2050032, 2020.
			

			
			\bibitem{article4} J. Bhowmick, D. Goswami and G. Landi: Levi-Civita connections and vector fields for noncommutative differential calculi, {Internat. J. Math.}, {\bf 31}, No 8, 2020.
     

					\bibitem{article6} J. Bhowmick and S. Mukhopadhyay: Covariant connections on bicovariant differential calculus, {J. Algebra}, {\bf 563}, 2020, 198--250.
      
					
					\bibitem{article7} J. Bhowmick, D. Goswami and S. Joardar: Levi-Civita connections for conformally deformed metrics on tame differential calculi, arXiv: 2101. 07221. 
					
		\bibitem{sitarz}	A. Bochniak, A. Sitarz and P. Zalecki: Riemannian geometry of a discretized circle and torus, arXiv:2007.01241.		
					
	\bibitem{chak_sinha} P. S. Chakraborty and K.B. Sinha: Geometry on the Quantum Heisenberg Manifold, Journal of Functional Analysis, {\bf 203}, 425--452, 2003.


 \bibitem{connes} A. Connes: Noncommutative geometry. Academic Press, San Diego, CA, 1994.
			
			
			\bibitem{Connes-dubois} A. Connes and M. Dubois-Violette : Noncommutative finite-dimensional manifolds. I.,
Spherical manifolds and related examples. Comm.\ Math.\ Phys. \textbf{230} (2002), no. 3, 539--579.





\bibitem{connes_landi} A. Connes and G. Landi: Noncommutative Manifolds the Instanton Algebra and Isospectral Deformations, Commun. Math. Phys.221:141--159,2001.


\bibitem{Connes_moscovici} A. Connes and H. Moscovici: Modular curvature for noncommutative two-tori, 
J. Amer. Math. Soc. {\bf 27 }(2014) 639--684.

\bibitem{scalar_3} A. Connes and P. Tretkoff: The Gauss-Bonnet theorem for the noncommutative two
torus, in Noncommutative Geometry, Arithmetic, and Related Topics, Johns Hopkins Univ. Press, Baltimore, MD, 2011, 141--158.

\bibitem{dubois} M. Dubois-Violette and P.W. Michor: Derivation et calcul differentiel non commutatif II, C. R. Acad. Sci. Paris Ser. I Math, {\bf 319} (1994) 927--931. 

\bibitem{dubois2} M. Dubois-Violette and P.W. Michor: Connections on central bimodules, J. Geom. Phys. {\bf 20} (1996) 218--232.


\bibitem{khalkhali} F. Fathizadeh, M. Khalkhali: Curvature in Noncommutative Geometry, arXiv:1901.07438.

 
\bibitem{frolich} J. Frohlich, O. Grandjean and A. Recknagel: Supersymmetric Quantum
Theory and Non-Commutative Geometry, Commun. Math. Phys {\bf 203} (1999) 119--184.

	
	\bibitem{heckenberger_etal}   I. Heckenberger and K. Schmuedgen : Levi-Civita Connections on the Quantum Groups 
	${\rm SL}_q(N)$ ${\rm O}_q(N)$ and ${\rm Sp}_q(N)$,
Comm. Math. Phys. {\bf 185} (1997) 177--196.

\bibitem{soumalya} S. Joardar: Scalar Curvature of a Levi-Civita Connection on Cuntz algebra with three generators, Lett. Math. Phys., {\bf 109}, 2665--2679, 2019. 



\bibitem{majid_1} S. Majid: Noncommutative Riemannian and spin geometry of the standard q-sphere, Commun.
Math. Phys. {\bf 256 }(2005) 255--285.

\bibitem{matassa} M. Matassa: Fubini-Study metrics and Levi-Civita connections on quantum projective spaces, arXiv: 2010. 03291. 

\bibitem{sheu} M.A. Peterka and A.J.L. Sheu:  On Noncommutative Levi-Civita Connections, International Journal of Geometric Methods in Modern Physics {\bf 14}, No. 5 (2017) 1750071.

\bibitem{rieffel_heisenberg}  M.A. Rieffel: Deformation quantization of Heisenberg manifolds, \emph{Comm.~ Math.~ Phys.} \textbf{122} (1989),  531--562.

\bibitem{rieffel}  M. A. Rieffel : Deformation Quantization for actions of $ \IR^{d},$ Memoirs of the American Mathematical Society, November 1993. Volume {\bf 106}, Number 506.

\bibitem{Rosenberg} J. Rosenberg: Levi-Civita's Theorem for Noncommutative Tori, SIGMA {\bf 9 }(2013),
071.

\bibitem{weber} T. Weber: Braided Cartan Calculi and Submanifold Algebras, arXiv: 1907.13609.


\bibitem{Skd}  M. Skeide: Hilbert modules in quantum electrodynamics and quantum probability.  Comm. Math. Phys. \textbf{192 } (1998), 569--604.


\end{thebibliography}
\end{document}